\documentclass[journal]{new-aiaa}
\usepackage[utf8]{inputenc}
\usepackage{textcomp}

\usepackage{graphicx}
\usepackage{amsmath}
\usepackage[version=4]{mhchem}
\usepackage{siunitx}
\usepackage{longtable,tabularx}
\usepackage{booktabs}

\usepackage{amsthm}
\usepackage{multirow}
\usepackage{subcaption}
\usepackage{rotating}

\newcommand{\sat}{\mathrm{sat}}
\newcommand{\wait}{\mathrm{wait}}
\newcommand{\pln}{\mathrm{plane}}
\newcommand{\prk}{\mathrm{parking}}
\newcommand{\hld}{\mathrm{hold}}
\newcommand{\lau}{\mathrm{lau}}
\newcommand{\inc}{i}
\newcommand{\mnv}{\mathrm{maneuv}}
\newcommand{\mnf}{\mathrm{manufac}}
\newcommand{\oos}{\mathrm{oos}}

\newcommand{\trs}{\mathrm{trans}}
\newcommand{\AMC}{\mathrm{AMC}}
\newcommand{\AP}{\mathrm{AP}}
\newcommand{\ofr}[1]{\rho_{#1}}
\newcommand{\stklv}[1]{\overline{\mathrm{SL}}_{#1}}
\newcommand{\estg}[1]{\mathrm{ES}_{#1}}
\newcommand{\acost}[1]{C_{#1}}
\newcommand{\scost}[1]{c_{#1}}
\newcommand{\ofq}[1]{n_{#1}}
\newcommand{\nbr}[1]{N_{#1}}
\newcommand{\alt}[1]{h_{#1}}
\newcommand{\mtm}[1]{T_{#1}}
\newcommand{\frate}[1]{\lambda_{#1}}
\newcommand{\rrate}[1]{\mu_{#1}}
\newcommand{\dec}[1]{\mathbf{x}_{#1}}
\newcommand{\prm}[1]{\mathbf{q}_{#1}}
\newcommand{\cdf}[1]{P_{#1}}
\newcommand{\pdf}[1]{f_{#1}}

\newcommand{\rev}[1]{{\color{black}{#1}}}

\title{On-Orbit Servicing-Integrated Maintenance Strategy for Satellite Constellation\footnotetext{Presented as Paper 2026-2631 at the AIAA SciTech Forum, Orlando, FL, January 12--16, 2026}}

\author{Jaewoo Kim \footnote{Ph.D. Candidate, Department of Aerospace Engineering, 291 Daehak-Ro. Student Member AIAA}}
\affil{Korea Advanced Institute of Science and Technology, Daejeon 34141, Republic of Korea}
\author{Taehyun Sung \footnote{Senior Researcher.}}
\affil{Agency for Defense Development, Daejeon 34186, Republic of Korea}
\author{Woonam Hwang \footnote{Associate Professor of Industrial and Systems Engineering, 291 Daehak-Ro.} and Jaemyung Ahn\footnote{Professor of Aerospace Engineering, 291 Daehak Ro; jaemyung.ahn@kaist.ac.kr. Associate Fellow AIAA (Corresponding Author).}}
\affil{Korea Advanced Institute of Science and Technology, Daejeon 34141, Republic of Korea}

\begin{document}

\maketitle

\begin{abstract}
This paper proposes a maintenance strategy for a satellite constellation that utilizes on-orbit servicing (OOS). Under this strategy, the constellation operator addresses satellite failures in two ways: by deploying new satellites and by recovering failed satellites through OOS. We develop an inventory management model with a parametric replenishment policy for the maintenance process, which can evaluate the performance of the satellite constellation system. \rev{Based on this model, we formulate the interaction between the constellation operator---who seeks to maintain the required service level of the constellation while minimizing maintenance cost---and the OOS provider---who seeks to maximize profit by selecting service price and performance levels---as a bi-objective optimization problem and identify the corresponding Pareto-optimal solutions. A case study based on real-world-scale constellation and launch service shows that, relative to the benchmark strategy without OOS, the OOS-integrated solutions can reduce annual maintenance cost by up to 14.5\%, while reducing annual launch and manufacturing costs by approximately 25\% each and maintaining the required service levels. The results further show that, within a given scenario, the Pareto-optimal set is generally generated by an almost invariant maintenance strategy on the constellation operator side, whereas most variation along the Pareto frontier is driven by pricing decisions of the OOS provider side. Among the OOS-related parameters, the fraction of failures that can be recovered through OOS has the strongest structural effect.}
\end{abstract}

\section*{Nomenclature}
{\renewcommand\arraystretch{1.0}
\noindent\begin{longtable*}{@{}l @{\quad=\quad} l@{}}
$\AMC$ & constellation operator's annual maintenance cost, \$M/year\\
$\AP$ & OOS provider's annual profits, \$M/year\\
$\acost{\hld}$ & annual holding cost, \$M/year\\
$\acost{\lau}$ & annual launch cost, \$M/year\\
$\acost{\mnv}$ & annual maneuvering cost, \$M/year\\
$\acost{\mnf}$ & annual manufacturing cost, \$M/year\\
$\acost{\oos}$ & annual OOS cost, \$M/year\\
$\scost{\hld}$ & annual holding cost of a satellite, \$M/satellite/year\\
$\scost{\lau}$ & launch costs, \$M\\
$\scost{\mnf}$ & manufacturing cost of a satellite, \$M/satellite\\
$\scost{\oos}$ & cost for the OOS provider to recover a satellite, \$M\\
$\estg{\prk}$ & expected shortage of parking spares, -\\
$\estg{\pln}$ & expected shortage of in-plane spares, in units of satellites\\
$\alt{\prk}$ & altitude of parking orbits, km\\
$\alt{\pln}$ & altitude of orbital planes in constellation, km\\
$\inc$ & inclination angle of an orbit, deg \\
$k_s$ & reorder point of parking spares, -\\
$k_Q$ & order quantity of parking spares, -\\
$M_{\mathrm{dry}}$ & dry mass of satellite, kg\\
$M_{\mathrm{fuel}}$ & fuel mass required for orbital transfer, kg\\
$\dot{M}_{\mathrm{prop}}$ & mass flow rate of satellite propulsion system, kg/s\\
$\nbr{\oos}$ & the maximum number of service for a satellite, -\\
$\nbr{\prk}$ & number of parking orbits, in units of planes\\
$\nbr{\pln}$ & number of orbital planes in constellation, in units of planes\\
$\nbr{\sat}$ & number of satellites per orbital plane in constellation, in units of satellites\\
$\nbr{t}$ & number of time units per year, in units of time per year\\
$\ofq{\oos}$ & annual service frequency, per year\\ 
$\ofq{\prk}$ & annual order frequency of parking spares, per year\\ 
$\ofq{\pln}$ & annual order frequency of in-plane spares, per year\\
$p_{\oos}$ & OOS price, \$M\\
$Q$ & order quantity of in-plane spares, in units of satellites \\
$Q_{\max}$ & launch capacity, in units of satellites\\
$\stklv{\prk}$ & mean stock level of parking spares, -\\
$\stklv{\pln}$ & mean stock level of in-plane spares, in units of satellites\\
$\stklv{\mathrm{wait}}$ & mean stock level of failed satellites waiting for service, in units of satellites\\
$r_{\oos}$ & fraction of serviceable in-orbit failures, - (or \%)\\
$s$ & reorder point of in-plane spares, in units of satellites \\
$\mtm{d}$ & mean time to disposal of satellites serviced up to the maximum, in time units\\
$\mtm{\mathrm{life}}$ & maximum lifespan of satellites, years \\
$\mtm{\oos}$ & mean time spent during awaiting OOS until disposal, in time units\\
$\mtm{\mathrm{opr}}$ & mean time spent during operations until disposal, in time units\\
$\mtm{\prk}$ & mean time spent as parking spares until disposal, in time units\\
$\mtm{\pln}$ & mean time spent as in-plane spares until disposal, in time units\\
$\mtm{\mathrm{trans}}$ & mean time spent during orbital transfer from parking to target orbit, in time units\\
$\mathrm{TOF}$ & time of flight, in time units\\
$t_{\trs}$ & in-plane transfer maneuver time, in units of time\\
$\mathcal{T}_{j}$ & admissible time interval for supplementation from the \rev{$j$th} closest parking orbit\\
$V_{\mathrm{e}}$ & exhaust velocity of satellite propulsion system, km/s\\
\rev{$\dec{}$} & \rev{decision variables vector of optimization problem}\\
$\alpha_1,\alpha_2$ & shape parameters of cost-responsiveness function, -\\
$\gamma_m$ & fraction of in-plane spares serviced $m$ times, - (or \%)\\ 
$\delta$ & random variable for stochastic demand \\
$\epsilon_{\mathrm{fuel}}$ & fuel mass conversion coefficient, \$M/kg\\ 
$\frate{\prk}$ & demand rate at parking spares, per unit time \\
$\frate{\prk}$ & demand rate at in-plane spares, in units of satellites per unit time \\
$\frate{\sat}$ & failure rate of a satellite, in units of satellites per year\\
$\rrate{\lau}$ & parameter of stochastic launch delay, in units of inverse time\\
$\rrate{\oos}$ & responsiveness of service, in units of inverse time\\ 
$\ofr{\prk}$ & order fill rate for parking spares, - (or \%)\\
$\ofr{\pln}$ & order fill rate for in-plane spares, - (or \%)\\ 
$\tau$ & random variable for stochastic lead time\\
$\mathbb{R}_{++}$ & positive real number set \\
$\mathbb{Z}_{+}$, $\mathbb{Z}_{++}$ & non-negative, and positive integer sets\\
\end{longtable*}}

\section{Introduction}\label{section 1}
\lettrine{S}{atellite} mega-constellations have played a pivotal role in shaping the modern space industry, driven by the demand for global coverage services and advances in space technologies. These driving forces have encouraged active participation from major players, including SpaceX's Starlink, \rev{Amazon Leo}, and Eutelsat's OneWeb \cite{pachler2021updated, pachler2024flooding}. This participation has expanded the market and established mega-constellations as a critical component in the ecosystem of the space industry.

\rev{From an operational perspective, mega-constellations in low Earth orbit (LEO) pose challenges throughout the lifecycle that differ substantially from those of traditional constellations in geostationary equatorial orbit (GEO). The large number of satellites, together with the relatively small size of individual spacecraft and tight constraints on power, performance, and maneuvering capability, complicates system operation and management. Prior research has addressed multiple lifecycle phases, including pattern design \cite{walker1971circular,ballard1980rosette,mortari2004flower,ullock1963optimum,beste1978design,hanson1989improved,casanova2014seeking,ma2020hybrid,huang2021multi,jia2022design,wang2022mega,xiao2025low,rogers2026optimal}, deployment \cite{cornara1999satellite,deweck2004staged,crisp2015launch,anderson2022design,pasquale2022optimization,sung2023optimal}, mission operations \cite{chen2019distributed,valicka2019mixed,kim2026integrated,mao2026spatiotemporal}, and end-of-life management \cite{huang2019large,olivieri2020large}. This study focuses on maintenance. In mega-constellations, component-level failures (i.e., in-orbit failures of individual satellites) may occur more frequently, not only because more satellites are in operation, but also because individual satellites are often designed under tighter mass, cost, and redundancy constraints. In addition, rapid technology refresh cycles may require agile integration of payload and spacecraft upgrades. For example, Starlink satellites are reported to have an operational lifetime of about five years, a wet mass of several hundred kilograms, and a planned constellation size of 42,000 \cite{pultarova2025starlink}. In contrast, current Intelsat---a GEO-based telecommunication company---satellites can exceed 6 metric tons, are expected to remain operational for more than 15 years \cite{krebs}, and are deployed in fleets of only a few tens of satellites. These differences motivate the development of a systematic maintenance strategy tailored to mega-constellations, specifically a structured \textit{logistics} framework for effective failure response \cite{ho2024space}. Recent studies, including Jakob et al. \cite{jakob2019optimal} and Kim et al. \cite{kim2025replenishment}, have proposed maintenance strategies and corresponding assessment tools for satellite mega-constellations.}

Specifically, Jakob et al. \cite{jakob2019optimal} proposed a multi-echelon supply chain for spare satellites, along with an accompanying inventory management policy and analytical forms of performance metrics. The supply chain comprises three echelons---ground spares, parking spares, and in-orbit spares---which are analogous to suppliers, warehouses, and retailers, respectively, in a terrestrial supply chain. Within this framework, parking spares and in-plane spares are replenished using distinct $\left(s,Q\right)$ policies \cite{silver2016inventory}, where a batch of $Q$ units is ordered once the stock level drops to $s$. Under this structure, the authors assumed a Poisson distribution for satellite failures and an exponential distribution for launch lead times. They derived analytical expressions for inventory-related parameters, maintenance costs, and service levels. The cost model encompasses several elements, including satellite manufacturing, launch, holding, and maneuvering. Finally, they formulated an optimization problem from the constellation operator's perspective to determine the policy parameters and parking orbit geometries that minimize maintenance costs while satisfying the required service level of the maintenance system.

Building on the work of Jakob et al., Kim et al. \cite{kim2025replenishment} established a dual supply strategy by introducing an auxiliary (direct) channel that directly injects spare satellites into the target orbits \cite{sung2023optimal}. They integrated an inventory policy denoted by $\left(R_1,R_2,Q_1,Q_2\right)$ with a time window for in-plane spares. Under this policy, when the stock level falls to $R_1$, $Q_1$ units are ordered from the parking spares via the primary (indirect) supply channel. If the stock level subsequently drops to $R_2$ within the predefined time window, the auxiliary (direct) supply channel is activated, and $Q_2$ units are ordered. They derived analytical expressions for inventory-related parameters and formulated two optimization problems from different perspectives: the satellite constellation operator and the auxiliary launch service provider (i.e., the transportation service provider for the direct channel). For a comprehensive review of maintenance for small satellite constellations, readers are referred to Jakob et al. \cite{jakob2019optimal} and Kim et al. \cite{kim2025replenishment}.

\rev{Beyond replacing failed satellites with new spare units, on-orbit servicing (OOS) can restore operational capability and is beginning to be commercialized \cite{li2019orbit, arney2021orbit, davis2019orbit, cavaciuti2022space}. Historically, OOS has focused on high-value space systems such as the International Space Station, the Hubble Space Telescope, and communication satellites in GEO, for which life extension can provide substantial economic value. Numerous studies have assessed the benefits of OOS for single large systems \cite{saleh2002space, lamassoure2002space, joppin2006orbit, long2007orbit, yao2013orbit, liu2021economic, kim2025system}, and operational missions such as MEV-1 \cite{henry2020intelsat} indicate its practical relevance.

More recently, studies on OOS for servicing multiple satellites have examined technology and architecture alternatives \cite{dujonchay2017quantification,sears2018impact,ho2020semi,dujonchay2022on,luu2022on,cunningham2026economic,falcone2026single}, as well as the optimization of service planning and system design \cite{zhou2015mission,bang2018two,verstraete2018geosynchronous,bang2019multitarget,hudson2020versatile,dujonchay2021framework,han2022multiple,han2023optimal,lee2023optimal,yan2025rapid,patnala2026on,zhu2020orbit,shimane2024orbital,choi2026orbital}. However, this body of literature primarily adopts the perspective of the OOS provider, focusing on how OOS systems---an emerging and still difficult-to-realize concept---should be designed and operated. By contrast, the customer-side perspective, namely how constellation operators should incorporate OOS into their system, remains relatively underexplored. For example, Luu and Hastings \cite{luu2022on} compared two maintenance strategies for satellite mega-constellations from the constellation operator's perspective: deploying spare satellites or relying on OOS. Nevertheless, to our knowledge, prior research has not jointly considered spare satellites and OOS within a unified maintenance strategy for satellite mega-constellations.

To address this gap, this paper proposes a maintenance strategy for satellite constellations that integrates OOS with spare satellite management. Under this strategy, failed satellites may be replaced by spare units or recovered through OOS when the failure is recoverable. To evaluate the proposed strategy, we develop an assessment tool based on an inventory management model. Building on this model, we define a decision-making scenario involving two stakeholders: the constellation operator, as the customer of OOS, and the OOS provider. We then formulate these settings as optimization problems and conduct a case study.

The contributions of this work are twofold. First, we propose an inventory control-based framework for integrating OOS into mega-constellation maintenance systems. The framework jointly considers spare-based replacement and recovery through OOS, and incorporates the resulting closed-loop structure into the inventory management model. To our knowledge, this closed-loop feature has not been addressed in prior work \cite{jakob2019optimal,kim2025replenishment}. Second, building on this framework, we formulate an optimization problem and present a case study to demonstrate its applicability and derive managerial insights for relevant decision makers. In particular, we conduct a parametric analysis of OOS cost and performance, both of which remain uncertain as OOS capabilities continue to emerge. From the customer perspective (i.e., constellation operators), the analysis guides how operators should respond to the emergence of OOS in the market. From the customer perspective (i.e., constellation operators), the analysis clarifies how operators should respond to the emergence of OOS in the market. From the provider perspective (i.e., OOS providers), it identifies the key factors that make OOS offerings attractive to customers and examines how customers may respond to such offerings.

The remainder of this paper is organized as follows. Section \ref{section 2} outlines the modeling assumptions and develops the supply chain and inventory management framework to evaluate the OOS-integrated maintenance strategy. Section \ref{section 3} validates the proposed model by comparing its results with those obtained from inventory simulations and assesses its robustness to key assumptions through analyses under relaxed assumptions. Section \ref{section 4} presents the decision-making scenarios and corresponding optimization problem formulations, and Section \ref{section 5} provides a case study illustrating their application. Finally, Section \ref{section 6} concludes the paper.}

\section{Model Construction}\label{section 2}
This section presents the proposed maintenance strategy and evaluation framework. \rev{We first describe the concept of operations of the maintenance system, together with the associated supply chain and its underlying assumptions. We then present the parametric models to characterize the mean inventory behavior at each echelon of the supply chain and relate these models to the maintenance cost and other system performance parameters.

\subsection{Concept of Operations}\label{section: 2-A}
Before presenting the detailed assumptions and inventory models, we briefly describe the maintenance system. Figure \ref{fig: conops} provides an overview of the proposed maintenance concept.

\begin{figure}[hbt!]
    \centering
    \includegraphics[page=1, width=.7\textwidth]{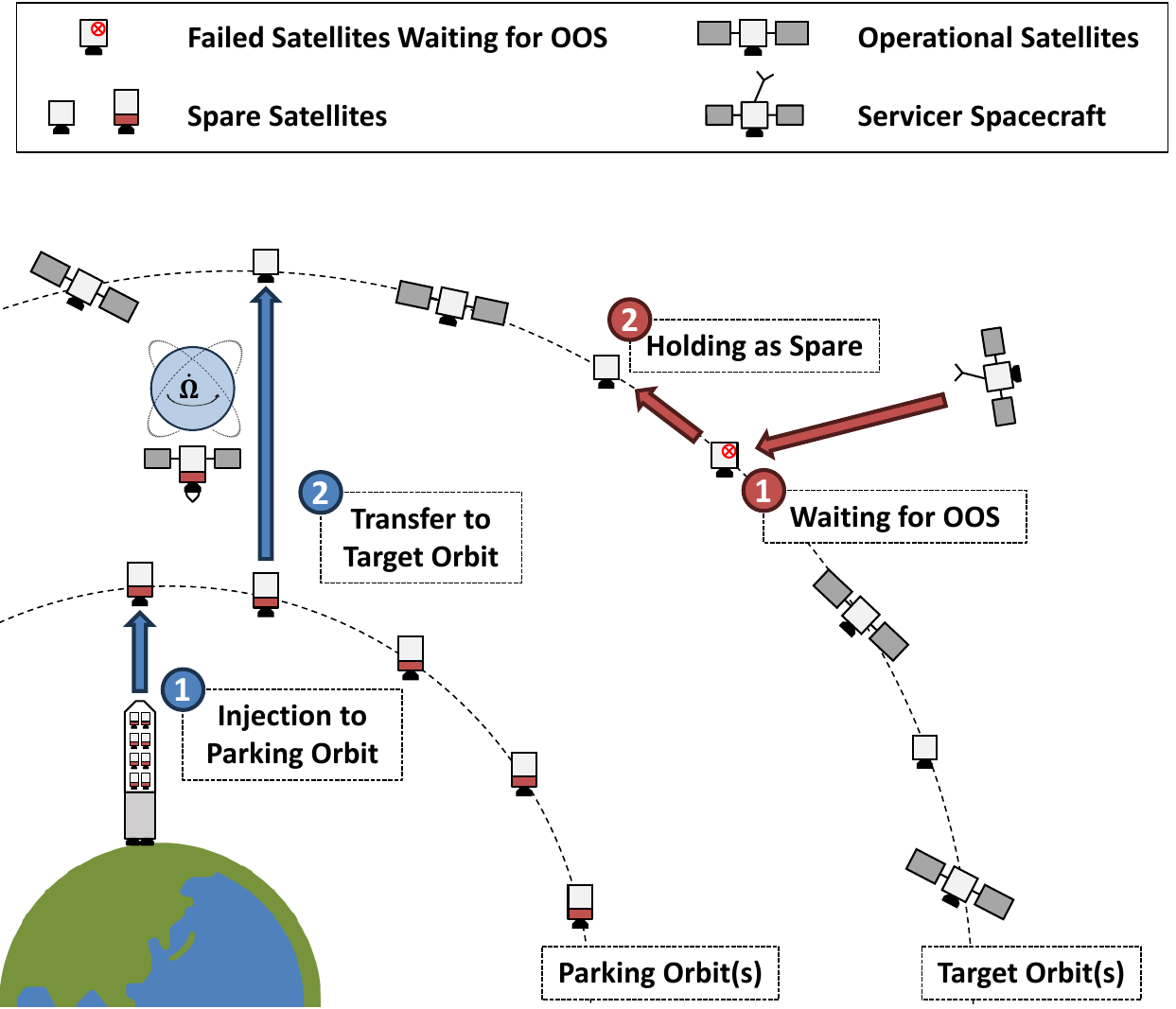}
    \caption{\rev{Illustration of the proposed OOS-integrated maintenance strategy.}}
    \label{fig: conops}
\end{figure}

We assume that the target constellation follows a Walker Delta pattern \cite{wertz2001mission}, where all orbital planes are circular, share the same inclination, and are uniformly distributed over the right ascensions of the ascending nodes (RAANs). Satellites within each orbital plane are also assumed to be uniformly distributed along, and all satellites in the constellation are assumed to be homogeneous. Although the formulation is developed for this symmetric configuration, it can be adapted with minimal effort to constellations that deviate from it, for example, by allowing nonzero eccentricity in each orbital plane or nonuniform distribution of orbital planes over the RAAN domain. In such cases, the model can still be used to track the mean behavior of the spare inventory, although deviations from the mean would require additional study.

To address failures of operational satellites, we consider three echelons of spare satellites, following Jakob et al. \cite{jakob2019optimal}. The first echelon consists of in-plane spares, which are distributed across the (near-)target orbits (i.e., the operational planes of the constellation). When an operational satellite fails, its slot is filled by an in-plane spare. The second echelon consists of parking spares, which replenish the in-plane spares. These parking orbits have the same inclination as the target orbits and are uniformly distributed over RAANs, but they are located at lower altitudes. A parking spare replenishes an in-plane spare by performing an orbit-raising maneuver once the corresponding parking and target orbits become aligned in RAAN due to the differential RAAN drift induced by the altitude difference, as discussed later. The third echelon consists of ground spares, which replenish the parking spares. Ground spares are newly manufactured satellites that are delivered to the parking orbits by launch vehicles. The associated replenishment flows of newly manufactured satellites are shown as blue arrows in Fig. \ref{fig: conops}.

In addition to the replenishment flow of newly manufactured satellites, the maintenance system includes an OOS-based recovery process, represented by the red arrows in Fig. \ref{fig: conops}. If an operational satellite experiences a failure mode that is recoverable through OOS, the failed satellite remains near the target orbit of its original slot and waits for service. A servicer spacecraft deployed in advance by the OOS provider then performs OOS, after which the failed satellite is retained as an in-plane spare.

\subsection{Supply Chain Overview and Assumptions}\label{section: 2-B}
The satellite constellation considered in the proposed strategy is assumed to follow a Walker Delta pattern with inclination $\inc$ and altitude $\alt{\pln}$. The constellation consists of $\nbr{\pln}$ orbital planes, whose RAANs are equally spaced by $2\pi/\nbr{\pln}$. Each orbital plane contains $\nbr{\sat}$ operational satellites.

As noted earlier, the constellation operator sources spare satellites through two channels: (1) recovering the operational life of failed satellites via OOS and (2) launching newly manufactured satellites from the ground. Figure \ref{fig: supply chain overview} provides an overview of the proposed supply chain.

\begin{figure}[hbt!]
    \centering
    \includegraphics[page=2, width=1.0\textwidth]{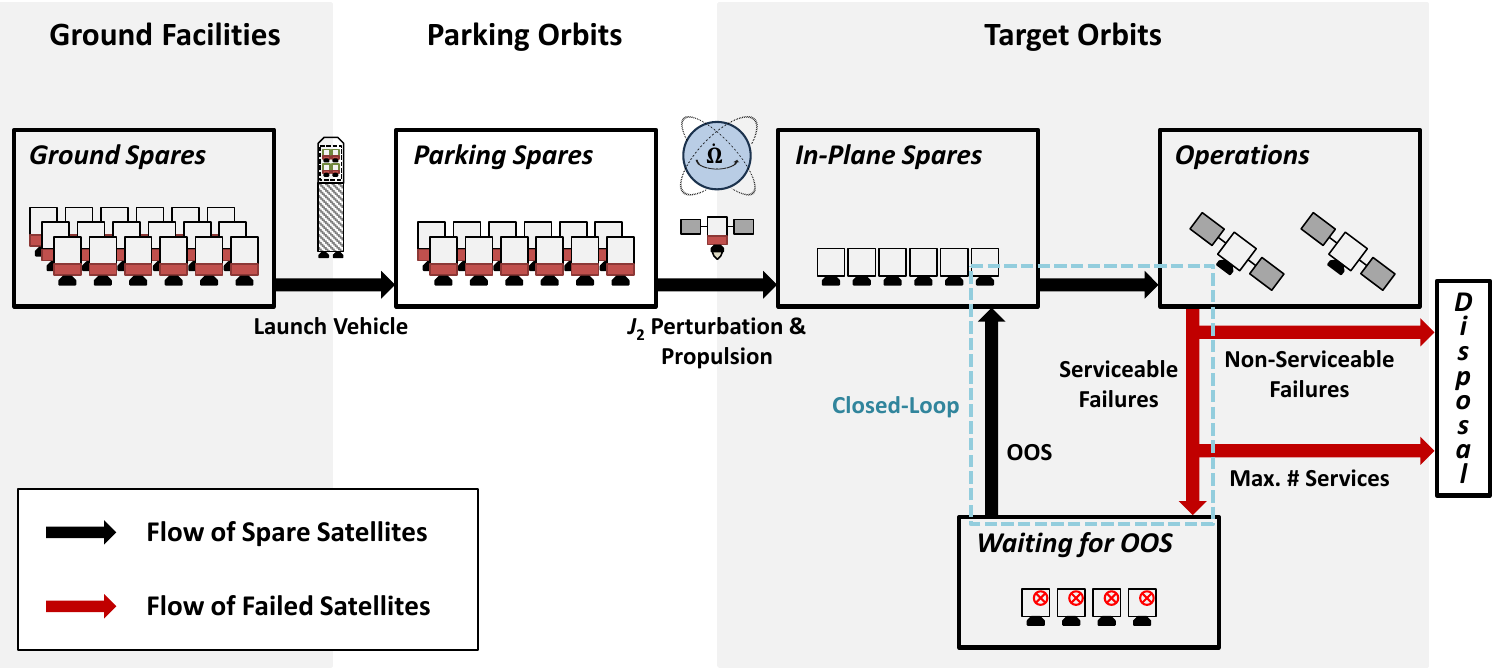}
    \caption{Overview of the supply chain for the OOS-integrated maintenance strategy.}
    \label{fig: supply chain overview}
\end{figure}
}
\subsubsection{Sourcing Method 1: Life Extension by OOS}
OOS addresses serviceable failures. The specific failure modes that are considered serviceable depend on the spacecraft design (e.g., service-friendly interfaces, degree of modularity), the maturity of OOS technologies, and the availability of supporting infrastructure. At the strategic level, we assume a fixed value for the fraction of serviceable in-orbit failures, $r_{\oos}$. In other words, if 100 failures occur, then on average $100r_{\oos}$ are serviceable and $100\left(1-r_{\oos}\right)$ are not. When a serviceable failure occurs, the failed satellite waits for servicing. After servicing, it is assumed to fully regain its operational capability and remain on standby as an in-plane spare.

In addition, even if a serviceable failure occurs, the operator may choose not to recover the satellite if other components, those that OOS cannot renew, are obsolete and require replacement. Therefore, we further assume that each satellite has an operational lifespan, denoted by $\mtm{\mathrm{life}}$. For example, if OOS can only refuel a depleted satellite, other components (e.g., electronics, structure) remain unchanged. In such cases, repeatedly refueling the satellite beyond its structural or technological lifespan may be inadequate. Consequently, during the strategic phase, the operator may define $\mtm{\mathrm{life}}$ based on the OOS capabilities and spacecraft design. As a proxy for $\mtm{\mathrm{life}}$, we introduce the maximum number of OOS per satellite, $\nbr{\oos}$, as a tracking parameter.

\rev{
In this study, we consider the aspects most relevant to the customer, namely, service price and responsiveness. Specifically, we assume that the OOS provider has already deployed the infrastructure required to support the constellation under consideration, including orbital depots for fuel, replacement units, and other servicing resources, as well as multiple servicing spacecraft capable of approaching and restoring failed satellites. Under this assumption, the OOS system has sufficient capacity to respond to service requests with stable cost and performance.

The OOS price is the amount the constellation operator pays to use the OOS. In this study, OOS is charged per servicing event at a fixed price of $p_{\oos}$. Correspondingly, the OOS provider receives revenue of $p_{\oos}$ per servicing event. To provide OOS, the provider must allocate resources to each servicing event, and we model the corresponding per-event cost as $\scost{\oos}$. As a result, the OOS provider's profit per servicing event is $\left(p_{\oos}-\scost{\oos}\right)$.

Service lead time is the representative measure of service responsiveness. It is defined as the duration between receipt of a service request and completion of the requested service, and it is the key performance parameter considered in this study. It consists primarily of the time required for the sequence of maneuvers by which a servicing spacecraft approaches the target satellite, including rendezvous, proximity operations, and docking, together with the time required to perform the servicing activity itself. The lead time depends on the current state of the OOS system. For example, if an available servicing spacecraft is far from the target satellite when the request is made, the resulting lead time is longer. Because the time required for a servicing spacecraft to reach the target is inherently uncertain, we model service lead time as a random variable. The modeling details are presented later in this section.

Note that the cost ($\scost{\oos}$) and the responsiveness of OOS are closely related. Given the OOS technologies, service types, and target failure modes of the serviced spacecraft, the OOS provider's operational effort determines both service responsiveness and cost. For example, if the OOS provider allocates more resources, such as additional servicing spacecraft, depots, and propellant for faster orbital transfers, the servicing cost increases, but so does responsiveness. This relationship is modeled in a simplified manner in this study and is introduced later in Section \ref{section 4}.
}

\subsubsection{Sourcing Method 2: Deploying New Satellites}
If a failure is non-serviceable or the failed satellite has exceeded its lifespan, a spare satellite replaces the failed one. They are managed by the strategy proposed by Jakob et al. \cite{jakob2019optimal}, which is the indirect channel in \cite{kim2025replenishment}. Spare satellites are distributed across three locations---ground facilities (ground spares), parking orbits (parking spares), and target orbits (in-plane spares). Here, the parking orbits share the same inclination as the target orbits, $\inc$, but have a lower altitude, $\alt{\prk}<\alt{\pln}$. There are $\nbr{\prk}$ parking planes, with RAANs equally spaced by $2\pi/\nbr{\prk}$.

\rev{In-plane spares directly replace failed satellites.} When the inventory of in-plane spares reaches a reorder point $s$, a batch of $Q$ spares is sourced from the parking spares, the next echelon. Here, the spare satellites are transported by leveraging $J_2$ perturbation to align RAANs of the parking and target orbits, and the satellite's propulsion system to raise the altitude. Here, we have the RAAN drift rate relation \cite{vallado2022fundamentals}:
\begin{equation}\label{eqn: nodal precession}
    \dot{\Omega}\left(a,e,i\right) =-\frac{3}{2} \sqrt{\frac{\mu_E}{a^3}} \frac{R_E^2}{a^2 (1-e^2)^2} J_2 \cos{i},
\end{equation}
where $a$, $\inc$, and $e$ are the semi-major axis, the inclination, and the eccentricity of the orbit, respectively, $J_2$ is the zonal harmonic coefficient, $R_E$ is the Earth's equatorial radius, and $\mu_E$ represents the Earth's gravitational constant.

Parking spares perform orbit-raising maneuvers to replenish the in-plane spares once the RAANs of the corresponding parking and target orbits are aligned. The motion is modeled using Keplerian dynamics under continuous thrust, reflecting the increasing adoption of low-thrust propulsion systems---such as electric thrusters---by small satellites due to their high efficiency. The required velocity increment ($\Delta V$) and time of flight ($\mathrm{TOF}$) are as follows \cite{alfano1994constant,bettinger2014mathematical}:
\begin{equation}\label{eqn: delta v}
    \Delta V = \sqrt{\frac{\mu_E}{r_{\mathrm{initial}}}} - \sqrt{\frac{\mu_E}{r_{\mathrm{final}}}},
\end{equation}
\begin{equation}\label{eqn: time of flight}
    \mathrm{TOF} = \frac{M_{\mathrm{fuel}}}{\dot{M}_{\mathrm{prop}}},
\end{equation}
where $r_{\mathrm{initial}}$ and $r_{\mathrm{final}}$ are the initial and final orbit radii, respectively ($r_{\mathrm{initial}} < r_{\mathrm{final}}$), $\mu_E$ is Earth's gravitational constant, $M_{\mathrm{fuel}}$ is the required fuel mass, and $\dot{M}_{\mathrm{prop}}$ is the mass flow rate of the propulsion system. Here, the fuel mass follows the rocket equation: 
\begin{equation}\label{eqn: rocket equation}
    M_{\mathrm{fuel}} = M_{\mathrm{dry}}(e^{\Delta V/V_{e}}-1),
\end{equation}
where $M_{\mathrm{dry}}$ is the satellite's dry mass and $V_{e}$ is the exhaust velocity.

Parking spares are managed in batches of size $Q$ to match the in-plane replenishment batch size. When the inventory of parking spares falls to $k_s$ batches, then $k_Q$ batches are newly launched from the ground to the parking orbit. Ground spares are assumed to be always available.

\rev{
We construct an inventory management model based on the assumptions listed below:
\begin{enumerate}
    \item In-orbit failures of operational satellites follow a Poisson process with constant rate $\frate{\sat}$. This assumption is reasonable under the steady-state operation of mega-constellations, where individual satellite health variations are averaged out at the fleet level, yielding an approximately constant aggregate failure rate.
    \item Failures of spare satellites are neglected; that is, parking and in-plane spares are assumed to be failure-free. This simplification is justified because the time spent in spare status is much shorter than the operational lifetime and can therefore be neglected at the strategic level.
    \item In-plane spares are assumed to be available instantaneously to replace a failed satellite. This is reasonable because the associated response time is short, typically on the order of a few days \cite{cornara1999satellite}, relative to other relevant time scales such as RAAN alignment, launch lead time, orbital transfer, and OOS lead time.
    \item Each target orbit is replenished from the closest available parking orbit, defined as the nonempty parking orbit with the earliest RAAN alignment. When the inventory of in-plane spare reaches the reorder point, $s$, a batch of size $Q$ is dispatched from that parking orbit.
    \item Ground-to-parking delivery is assumed to supply only a single parking orbit. Without an orbital transfer vehicle or additional mechanisms, supplying multiple parking orbits would be inefficient \cite{lang1998comparison}.
    \item A new in-plane replenishment order is placed only after the previous order has completed transit to facilitate the tracking of the orders.
    \item Stock-outs are assumed to be rare. Accordingly, the probability that all parking orbits are simultaneously out of stock is treated as negligible, consistent with nominal maintenance system operations.
    \item The reorder point and order quantity at parking orbits are restricted to integer multiples of the batch size ($Q$).
    \item The OOS system is assumed to have sufficient capacity to serve the target constellation, so that individual service requests may be treated independently and significant waiting time is rare. Under this assumption, OOS lead time is modeled as an exponential random variable. This choice is motivated by its common use for lead times with these characteristics and by its mathematical tractability \cite{silver2016inventory,ebeling2019introduction}. The corresponding mean time to recovery, or equivalently mean time to repair (MTTR), is given by $\rrate{\oos}^{-1}$.
\end{enumerate}

Assumptions 1--8 are adapted from prior works \cite{jakob2019optimal,kim2025replenishment}. These assumptions are appropriate for strategic-level analysis of mega-constellations because they preserve the key operational characteristics while reducing the model's complexity. Assumption 9 is newly introduced in this study. In the model validation (Section \ref{section 3}), we relax this assumption to examine whether the proposed model remains robust under alternative lead time conditions.
}

\subsection{In-Plane Spares}
In-plane spares directly respond to failures of operational satellites, which follow a Poisson process with rate parameter $\frate{\sat}$. The failure rate within a constellation's operational plane is defined as
\begin{equation}\label{eqn: failure rate in-plane}
    \frate{\pln}=\frac{\frate{\sat}\nbr{\sat}}{\nbr{t}},
\end{equation}
where $N_{\sat}$ is the number of satellites per plane, and $N_t$ is the number of discretized time units per year \cite{jakob2019optimal}.

\begin{figure}[hbt!]
    \centering
    \includegraphics[page=3, width=1.0\textwidth]{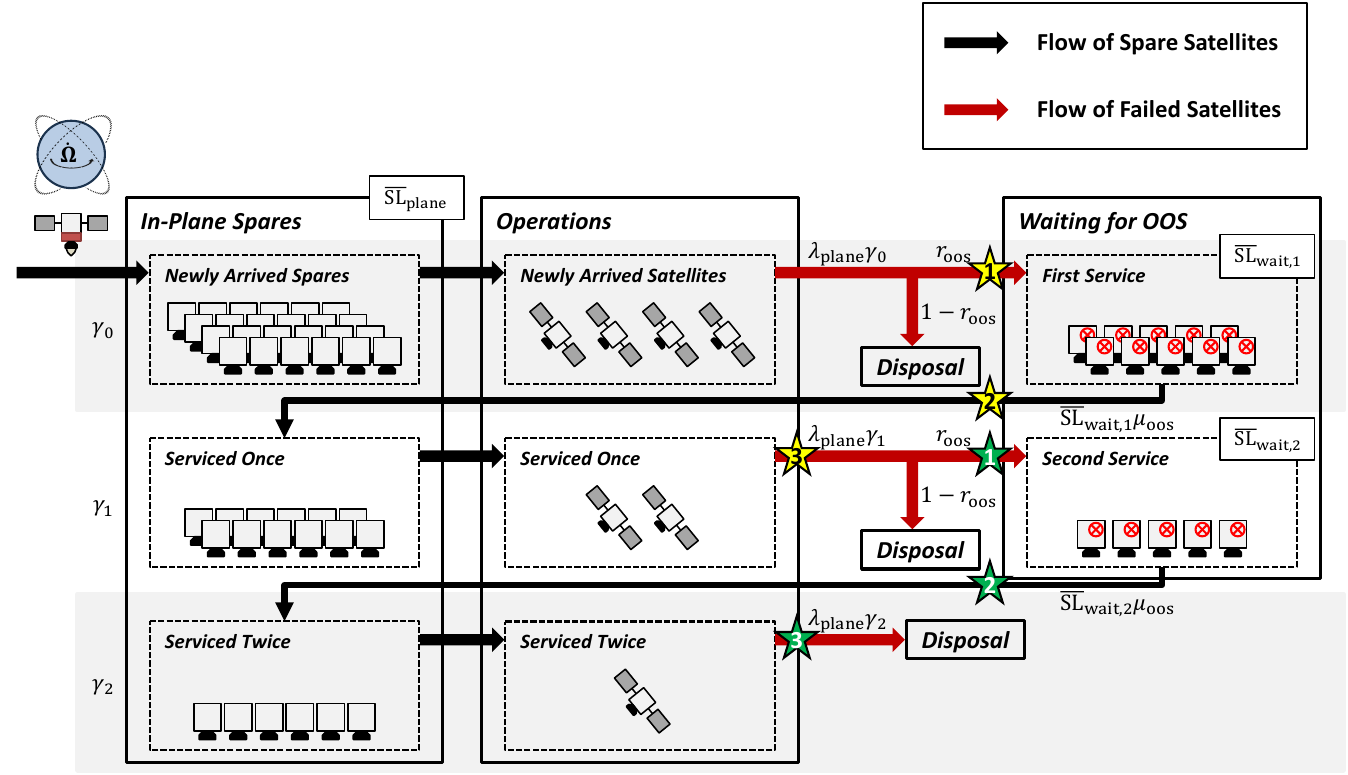}
    \caption{\rev{Example spare/failed satellites flows within a target orbit, $\nbr{\oos}=2$.}}
    \label{fig: example flows}
\end{figure}

In-plane spares are replenished by parking spares and OOS. Figure \ref{fig: example flows} provides a closer view of the closed-loop system shown in Fig. \ref{fig: supply chain overview}, illustrating the flows of spare and failed satellites for the case of $\nbr{\oos} = 2$\rev{, meaning that a satellite can be recovered at most twice. From the perspective of replacing the slot vacated by a failed satellite, a spare newly transferred from the parking orbit and a spare recovered through OOS are operationally equivalent. However, for modeling purposes, we distinguish in-plane spares according to their service history, namely whether they are newly arrived from the parking orbit or have been recovered through OOS, and, in the latter case, by the number of times they have been serviced.

After a parking spare arrives in the target orbit, it remains as an in-plane spare until an operational satellite fails and vacates a slot. The spare then fills that slot and operates until it fails. If the failure is serviceable, which occurs for a fraction $r_\oos$ of failures, the failed satellite waits for OOS and returns as an in-plane spare after the service lead time. Otherwise, it is disposed of. This process repeats until the satellite experiences a non-serviceable failure or reaches the maximum allowable number of services, $\nbr{\oos}$.

Over an infinite horizon, let $\gamma_0,\gamma_1,\dots,\gamma_{\nbr{\oos}}$ denote the steady-state fractions of in-plane spares. Here, $\gamma_0$ denotes newly arrived spares from parking orbits, while $\gamma_m$ ($m=1,2,\dots,\nbr{\oos}$) denotes spares that have undergone service exactly $m$ times. Because in-plane spares wait to replace failed operational satellites, they are held as spare inventory. We denote their mean stock level by $\stklv{\pln}$, which includes both newly arrived spares from parking orbits and spares that have been serviced any number of times. Likewise, serviceable failed satellites are held while awaiting service. We denote the corresponding mean stock level for satellites awaiting the $m$th service by $\stklv{\wait,m}$.
}

Assuming that the service lead time follows an exponential distribution with the mean value of $\rrate{\oos}^{-1}$ implies that the OOS process recovers a failed satellite with probability $\rrate{\oos}$ per unit time \cite{ebeling2019introduction}. \rev{The corresponding steady-state flow balance equations over an infinite horizon for 1) failed satellites awaiting the $m$th OOS and 2) in-plane spares that have been serviced $m$ times are indicated in Fig. \ref{fig: example flows} by the star pairs labeled 1-2 and 2-3, respectively. Within each color group (yellow and green), adjacent star pairs represent individual flow-balance relations, corresponding to Eqs. \eqref{eqn: flow balance 1} and \eqref{eqn: flow balance 2}:}
\begin{equation}\label{eqn: flow balance 1}
    \frate{\pln}\gamma_{m-1}r_{\oos}=\stklv{\wait,m}\rrate{\oos},\text{ for }m=1,2,\dots,\nbr{\oos},
\end{equation}
\begin{equation}\label{eqn: flow balance 2}
    \stklv{\wait,m}\rrate{\oos}=\frate{\pln}\gamma_m,\text{ for }m=1,2,\dots,\nbr{\oos},
\end{equation}
where $\stklv{\wait,m}$ is the mean stock level of failed satellites that have been serviced $m$ times and are awaiting the next service. Also, the fractions $\gamma_m$ satisfy the normalization condition:
\begin{equation}
    \sum_{m=0}^{\nbr{\oos}}\gamma_m=1.
\end{equation}
Solving these balance equations determine the fractions $\gamma_{m}$ ($m=0,1,\dots,\nbr{\oos}$) as follows:
\begin{equation}\label{eqn: gamma0}
    \gamma_0=\frac{1-r_{\oos}}{1-\left(r_{\oos}\right)^{\nbr{\oos}+1}},
\end{equation}
\begin{equation}
    \gamma_{m}=r_{\oos}\gamma_{m-1},\text{ for }m=1,2,\dots,\nbr{\oos}.
\end{equation}
In addition, we define $\stklv{\wait}$ as follows:
\begin{equation}
    \stklv{\wait}=\sum_{m=1}^{\nbr{\oos}}\stklv{\wait,m}.
\end{equation}

The demand and recovery processes are modeled as Poisson processes. Since the stock drop corresponds to demand minus recovery, the resulting stock drop process follows a Skellam distribution \cite{skellam1946frequency}. Accordingly, the probability mass function (PMF) of the stock drop over a time interval $\tau$ is given by:
\begin{equation}
    \cdf{\pln}\left(\delta;\tau\right)
    = e^{-\frate{\pln}\tau\left(2-\gamma_0\right)}
    \frac{1}{\left(1-\gamma_0\right)^{\delta/2}}
    I_{\delta}\left(2\frate{\pln}\tau\sqrt{1-\gamma_0}\right),
    \quad \delta \in \mathbb{Z}_+,
\end{equation}
where $I_{\delta}$ denotes the modified Bessel function of the first kind \cite{abramowitz1965handbook}.

To respond to this stock level drop, we adopt an $\left(s,Q\right)$ replenishment policy: whenever the in-plane spare inventory level drops to the reorder point $s$, an order of size $Q$ is placed. The lead time for such an order depends on the availability of parking spares, and the time for RAAN alignment between the target and parking orbits and for orbit raising \cite{jakob2019optimal, kim2025replenishment}. Denoting the lead time by $\tau$, its probability density function (PDF) is given as:
\begin{equation}
    \pdf{\pln}\left(\tau\right)=P'_{\mathrm{av}}\left(j\right)\frac{\nbr{\prk}\left|\dot{\Omega}_{\mathrm{rel}}\right|}{2\pi},\text{ for }\tau\in \mathcal{T}_j,\ j=1,2,\dots,N_{\prk},
\end{equation}
where $\cdf{\mathrm{av}}'\left(j\right)$ is the normalized probability that the $j$th closest parking orbit is the closest available parking one, $\dot{\Omega}_{\mathrm{rel}}$ is the relative RAAN drift rate between target and parking orbits, and \rev{$\mathcal{T}_j$ denotes the admissible range of lead times for sourcing from the $j$th closest parking orbit, including both the alignment time between the target orbit and its $j$th closest parking orbit and the orbit-raising time from the parking orbit to the target orbit, defined as}:
\begin{equation}
    \cdf{\mathrm{av}}'\left(j\right)=\frac{\cdf{\mathrm{av}}\left(j\right)}{\sum_{j'=1}^{\nbr{\prk}}\cdf{\mathrm{av}}\left(j'\right)},\quad \cdf{\mathrm{av}}\left(j\right)=\left(1-\ofr{\prk}\right)^{j-1}  \ofr{\prk},\text{ for }j=1,2,\dots,\nbr{\prk},
\end{equation}
\begin{equation}
    \dot{\Omega}_{\mathrm{rel}}=\dot{\Omega}\left(\alt{\prk},0,i\right)-\dot{\Omega}\left(\alt{\pln},0,i\right),
\end{equation}
\begin{equation}
    \mathcal{T}_{j}=\left[\left(j-1\right) \frac{2\pi}{\alt{\prk}} \frac{1}{\left|\dot{\Omega}_{\mathrm{rel}}\right|}+t_{\trs},j \frac{2\pi}{\nbr{\prk}} \frac{1}{\left|\dot{\Omega}_{\mathrm{rel}}\right|}+t_{\trs}\right),\text{ for } j=1,2,\dots,\nbr{\prk}.
\end{equation}
Here, $t_{\trs}$ is the time of flight from a parking orbit to a target orbit calculated by Eq. \eqref{eqn: time of flight}, \rev{and the right-open interval of $\mathcal{T}_j$ is used to avoid overlap between adjacent admissible intervals, thereby preventing double counting at the boundary.}

Following these demand and supply models, we define several inventory parameters. The mean stock level ($\stklv{\pln}$) represents the average inventory held. The order frequency ($\ofq{\pln}$) denotes how often orders are placed. The expected shortage ($\estg{\pln}$) is the average shortage within a reorder cycle, and the order fill rate ($\ofr{\pln}$) is the complement of the ratio between expected shortage and order quantity, representing the service level of the maintenance system. The service frequency ($\ofq{\oos}$) indicates the frequency of service calls. These are expressed as follows \cite{jakob2019optimal, kim2025replenishment, silver2016inventory}:
\begin{equation}
    \stklv{\pln}=\int_{\tau\ge t_{\trs}}\left(s-\frate{\pln}\gamma_0\tau+\frac{Q}{2}+\frac{1}{2}\right)\pdf{\pln}\left(\tau\right)\,d\tau,
\end{equation}
\begin{equation}
    \ofq{\pln}=\frac{\frate{\pln}\gamma_0}{Q}\nbr{t},
\end{equation}
\begin{equation}
    \estg{\pln}=\int_{\tau\ge t_{\trs}}\sum_{\delta\ge s}\left(\delta-s\right)\cdf{\pln}\left(\delta;\tau\right)\pdf{\pln}\left(\tau\right)\,d\tau,
\end{equation}
\begin{equation}
    \ofr{\pln}=1-\frac{\estg{\pln}}{Q},
\end{equation}
\begin{equation}
    \ofq{\oos}=\stklv{\wait}\rrate{\oos}\nbr{t}.
\end{equation}

\subsection{Parking Spares}
Parking spares respond to orders from in-plane spares, which are propagated from the failures of operational satellites. Empirically, if $\nbr{\pln}$ is large enough, a Poisson process can approximate the demand process in parking spares \cite{kim1991modeling,ganeshan1999managing}. Consequently, the demand rate for parking spares ($\frate{\prk}$) and the PMF of demand given time interval $\tau$ ($\cdf{\pln}$) are as follows:
\begin{equation}
    \frate{\prk}=\frac{\frate{\pln}\gamma_0\nbr{\pln}}{Q\nbr{\prk}},
\end{equation}
\begin{equation}
    \cdf{\prk}\left(\delta;\tau\right)=\frac{\left(\frate{\prk}\tau\right)^{\delta}}{\delta!}e^{-\frate{\prk}\tau},\quad \delta\in\mathbb{Z}_+.
\end{equation}
Similar to in-plane spares, parking spares are managed under an $(s,Q)$ policy with a reorder point $s=k_s$ and an order quantity $Q=k_Q$. Orders are placed in batches of size $k_Q$ and fulfilled from ground spares via launch vehicles. The lead time consists of a fixed processing time, $t_{\lau}$, and a stochastic delay due to launch scheduling, modeled as an exponential distribution with rate parameter $\rrate{\lau}$ \cite{jakob2019optimal}. As a result, the PDF of lead time is as follows:
\begin{equation}
    \pdf{\prk}\left(\tau\right)=\rrate{\lau}e^{-\left(\tau-t_{\lau}\right)\rrate{\lau}}.
\end{equation}

Following these demand and supply relationships, the key inventory metrics for parking spares, including mean stock level ($\stklv{\prk}$), order frequency ($\ofq{\prk}$), expected shortage ($\estg{\prk}$), and order fill rate ($\ofr{\prk}$), are computed as follows \cite{jakob2019optimal, kim2025replenishment, silver2016inventory}:
\begin{equation}
    \stklv{\prk}=\int_{\tau\ge t_{\lau}}\left(k_s-\frate{\prk}\tau+\frac{k_Q}{2}+\frac{1}{2}\right)\pdf{\prk}\left(\tau\right)\,d\tau,
\end{equation}
\begin{equation}
    \ofq{\prk}=\frac{\frate{\prk}}{k_Q}N_t,
\end{equation}
\begin{equation}
    \estg{\prk}=\int_{\tau\ge t_{\lau}}\sum_{\delta\ge k_s}\left(\delta-k_s\right)\cdf{\prk}\left(\delta;\tau\right)\pdf{\prk}\left(\tau\right)\,d\tau,
\end{equation}
\begin{equation}
    \ofr{\prk}=1-\frac{\estg{\prk}}{k_Q}\rev{,}
\end{equation}
\rev{where $N_t$ is the number of discretized time units per year.}

\subsection{Annual Maintenance Cost and Mean Time to Disposal}
Based on the inventory parameters for in‐plane and parking spares, we model the annual maintenance cost (AMC) for maintaining the satellite constellation. The AMC consists of the annual launch cost ($\acost{\lau}$), maneuvering cost ($\acost{\mnv}$), manufacturing cost ($\acost{\mnf}$), holding cost ($\acost{\hld}$)---which together represent the cost elements for managing spare satellites by previous works \cite{jakob2019optimal, kim2025replenishment}---and the servicing cost ($\acost{\oos}$):
\begin{equation}
    \AMC=\acost{\lau}+\acost{\mnv}+\acost{\mnf}+\acost{\hld}+\acost{\oos},
\end{equation}
where the individual cost components are as follows:
\begin{equation}
    \acost{\lau}=\scost{\lau}\ofq{\prk}\nbr{\prk},
\end{equation}
\begin{equation}
    \acost{\mnf}=M_{\mathrm{fuel}}\epsilon_{\mathrm{fuel}}\ofq{\pln}\nbr{\pln}Q,
\end{equation}
\begin{equation}
    \acost{\mnf}=\acost{\sat}\frate{\pln}\gamma_0\nbr{\pln}\nbr{t},
\end{equation}
\begin{equation}
    \acost{\hld}=\scost{\hld}\left(\stklv{\prk}Q\nbr{\prk}+\stklv{\pln}\nbr{\pln}+\stklv{\wait}\nbr{\pln}\right),
\end{equation}
\begin{equation}\label{eqn: oos cost}
    \acost{\oos}=p_{\oos}\ofq{\oos}\nbr{\pln}.
\end{equation}
Here, the annual launch cost is determined by the launch cost ($\scost{\lau}$) and the order frequency of parking spares ($\ofq{\prk}$); the annual maneuvering cost depends on the fuel mass required for orbit raising ($M_{\mathrm{fuel}}$), the specific cost of fuel ($\epsilon_{\mathrm{fuel}}$), and the order frequency of in-plane spares ($\ofq{\pln}$); the annual manufacturing cost is governed by the manufacturing cost per satellite ($\scost{\mnf}$), the failure rate ($\frate{\pln}$), and the fraction of in-plane spare demand met by newly launched satellites ($\gamma_0$); the annual holding cost depends on the holding cost rate ($\scost{\hld}$) and the mean stock levels at parking spares ($\stklv{\prk}$), in-plane spares ($\stklv{\pln}$), and failed satellites awaiting service ($\stklv{\wait}$); finally, the annual servicing cost is determined by the OOS price ($p_{\oos}$) and the service frequency ($\ofq{\oos}$).

Other than cost elements, the mean time from deployment to disposal of a satellite serviced $\nbr{\oos}$ times, $\mtm{d}$, is tracked. Suppose a satellite is disposed of after reaching the maximum number of allowable services. Then, it undergoes the following sequence: it remains as a parking spare after launch, transfers from a parking orbit to a target orbit, stays as an in-plane spare, operates, and is serviced after failure, repeatedly, until the service limit is reached.

Accordingly, $\mtm{d}$ can be modeled as the sum of the mean times spent in various phases: as a parking spare ($\mtm{\prk}$), during orbital transfer from parking to target orbit including RAAN alignment and orbit raising ($\mtm{\trs}$), as an in-plane spare ($\mtm{\pln}$), during operation ($\mtm{\mathrm{opr}}$), and while waiting for on-orbit servicing ($\mtm{\oos}$):
\begin{equation}
    \mtm{d}=\mtm{\prk} + \mtm{\trs} + \mtm{\pln} + \mtm{\mathrm{opr}} + \mtm{\oos},
\end{equation}
where
\begin{equation}
    \mtm{\prk} = \frac{\stklv{\prk}}{\frate{\prk}},
\end{equation}
\begin{equation}
    \mtm{\trs} = \int_{\tau\ge t_{\trs}}\tau \pdf{\pln}\left(\tau\right)\,d\tau,
\end{equation}
\begin{equation}
    \mtm{\pln} = \frac{\stklv{\pln}\left(\nbr{\oos}+1\right)}{\frate{\pln}},
\end{equation}
\begin{equation}
    \mtm{\mathrm{opr}} = \frac{\left(\nbr{\oos}+1\right)\nbr{t}}{\frate{\sat}},
\end{equation}
\begin{equation}
    \mtm{\oos}=\frac{\nbr{\oos}}{\rrate{\oos}}
\end{equation}
\rev{Note that $\mtm{\prk}$, $\mtm{\pln}$, and $\mtm{\oos}$ are obtained as the ratios of the corresponding stock levels to their outflow rates, $\mtm{\trs}$ follows from the definition of the mean of a continuous random variable with probability density function $f_{\pln}$, and $\mtm{\mathrm{opr}}$ is given by the reciprocal of the satellite failure rate. The factors ($\nbr{\oos}+1$) and $\nbr{\oos}$ arise because the satellite passes through the in-plane spare and operational phases ($\nbr{\oos}+1$) times, whereas it undergoes on-orbit servicing only $\nbr{\oos}$ times before disposal.}

\rev{
\section{Model Validation}\label{section 3}
The analytical model for the performance parameters constructed in Section \ref{section 2} is validated by comparing its values with those obtained from discrete simulations. We first describe the generation of the validation instances and then present the results.

\subsection{Instance Generation}
To evaluate the accuracy of the proposed inventory management model, we compared its estimates with the sample means obtained from $100$ Monte Carlo simulation runs. Each simulation spans $60$ years and is discretized into one-day intervals.

Some parameters were held fixed, while others were randomly sampled. Table \ref{tab: fixed parameters} lists the fixed parameters, and Table \ref{tab: sampled trade space} gives the ranges for the sampled parameters. Most parameter values were adopted from previous work \cite{kim2025replenishment}, whereas the OOS-related parameters were defined in this study to span a broad range of plausible values. The OOS price, $p_{\oos}$, was set to \$$2$M, which is approximately equal to the sum of the satellite production cost, $\scost{\sat}=$\$$0.5$M, and the per-satellite launch cost assuming a launch capacity of $40$ satellites, i.e., $\scost{\lau}/40\approx$\$$1.7$M. The same launch capacity was adopted as the reference launcher specification in the case study, following previous work \cite{kim2025replenishment}. The serviceable failure fraction, $r_{\oos}$, was varied from $0.1$ to $0.9$. The OOS MTTR, $\rrate{\oos}^{-1}$, was varied from $1$ week to $3$ months to represent a wide range of plausible service-performance levels in LEO \cite{luu2020valuation}. The maximum number of services, $\nbr{\oos}$, was varied from $1$ to $4$, based on the feasible range under the fixed satellite lifetime of $\mtm{\mathrm{life}}=30$ years.

\begin{table}[hbt!]
    \centering
    \caption{Fixed parameters for validation study}
    \begin{tabular}{lcc}
        \hline\hline
        Parameter, Unit & Notation & Value \\\hline
        Fuel mass conversion coefficient, \rev{\$M}/kg & $\epsilon_{\mathrm{fuel}}$ & 0.01 \\
        Satellite production cost, \rev{\$M}/satellite & $\scost{\sat}$ & 0.5\\
        Spare/failed satellite holding cost, \rev{\$M}/satellite/year & $\scost{\hld}$ & 0.5 \\
        Launch cost, \rev{\$M}/launch & $\scost{\lau}$ & 67 \\
        Specific impulse of satellite's propulsion system, sec & $I_{\mathrm{sp}}$ & 1200 \\
        Propellant mass rate, kg/s & $\dot{M}_{\mathrm{fuel}}$ & 1.3E-6 \\
        Satellite's dry mass, kg & $M_{\mathrm{dry}}$ & 150 \\
        OOS price, \rev{\$M} & $p_{\oos}$ & 2 \\
        Lifespan of satellite, years & $\mtm{\mathrm{life}}$ & 30 \\
        \hline\hline
    \end{tabular}
    \label{tab: fixed parameters}
\end{table}

\begin{table}[hbt!]
    \centering
    \caption{Sampled trade space of varying parameters for validation study}
    \begin{tabular}{lcc}
        \hline\hline
        Parameter, Unit & Notation & Range\\\hline
        Launch order processing time, weeks & $t_{\lau}$ & $\left[4, 16\right]$ \\
        Mean waiting time to launch, weeks & $\rrate{\lau}^{-1}$ & $\left[4, 16\right]$ \\
        Altitude of orbital planes in the constellation, km & $\alt{\pln}$ & $\left[500, 2000\right]$ \\
        Altitude of parking orbits, km & $\alt{\prk}$ & $\left[400, 1000\right]$ \\
        Inclination, deg & $\inc$ & $\left[40, 80\right]$\\
        Satellite failure rate, failures/year & $\frate{\sat}$ & $\left[0.01, 0.2\right]$ \\
        Number of orbital planes in the constellation, planes & \rev{$\nbr{\pln}$} & $\left[20, 40\right]$ \\
        Number of parking orbits, planes & $\nbr{\prk}$ & $\left[1, 20\right]$ \\
        Number of satellites per orbital plane in the constellation, satellites/plane & $\nbr{\sat}$ & $\left[20, 60\right]$ \\
        reorder point of in-plane spare, satellites & $s$ & $\left[1, 20\right]$ \\
        Order quantity at in-plane spare, satellites & $Q$ & $\left[1, 40\right]$ \\
        reorder point of parking spare, - & $k_s$ & $\left[1, 20\right]$ \\
        Order quantity at parking spare, - & $k_Q$ & $\left[1, 40\right]$ \\
        MTTR of OOS, weeks & $\rrate{\oos}^{-1}$ & $\left[1, 12\right]$ \\
        Serviceable failure fraction, - & $r_{\oos}$ & $\left[0.1, 0.9\right]$ \\
        Maximum number of services, - & $\nbr{\oos}$ & $\left[1,4\right]$ \\
        \hline\hline
    \end{tabular}
    \label{tab: sampled trade space}
\end{table}

 We repeatedly sampled the varying parameters until we obtained $100$ instances satisfying the following conditions:
\begin{equation}\label{eqn: val cond 1}
    \ofr{\pln}\ge \ofr{\pln}^{\mathrm{req}},
\end{equation}
\begin{equation}\label{eqn: val cond 2}
    \ofr{\prk}\ge \ofr{\prk}^{\mathrm{req}},
\end{equation}
\begin{equation}\label{eqn: val cond 3}
    \mtm{d}/\nbr{t}\le \mtm{\mathrm{life}},
\end{equation}
\begin{equation}\label{eqn: val cond 4}
    s \le Q,
\end{equation}
\begin{equation}\label{eqn: val cond 5}
    k_s \le k_Q.
\end{equation}
Here, Eqs. \eqref{eqn: val cond 1} and \eqref{eqn: val cond 2} represent the rare shortage assumption (Assumption 7) in Section \ref{section: 2-B} and are also used for service level requirements in later sections. Equation \eqref{eqn: val cond 3} imposes the lifespan condition, and Eqs. \eqref{eqn: val cond 4} and \eqref{eqn: val cond 5} maintain the cyclic behavior of the inventory. For $\ofr{\pln}^{\mathrm{req}}$ and $\ofr{\prk}^{\mathrm{req}}$, we considered values of $0.98$, $0.95$, and $0.90$, which yielded three sets of $100$ instances for testing stress cases under the rare shortage assumption.

\subsection{Results}
We tracked the mean stock levels of in-plane spares ($\stklv\pln$), parking spares ($\stklv\prk$), and failed satellites awaiting service ($\stklv{\wait}$); the order frequencies at the in-plane and parking spare echelons ($\ofq{\pln}$ and $\ofq{\prk}$); the OOS frequency ($\ofq{\oos}$); the mean time to disposal for satellites serviced up to the maximum allowable number of times ($\mtm{d}$); the annual maintenance cost ($\AMC$); and the order fill rates of the in-plane and parking spare echelons ($\ofr{\pln}$ and $\ofr{\prk}$).

For the rate parameters, $\ofr{\pln}$ and $\ofr{\prk}$, we used the absolute error:
\begin{equation}\label{eqn: model rate error}
\left|\mathrm{value}_{\mathrm{sim}} - \mathrm{value}_{\mathrm{model}}\right|,
\end{equation}
where $\mathrm{value}_{\mathrm{sim}}$ and $\mathrm{value}_{\mathrm{model}}$ are the parameter values obtained from the simulation and the model, respectively. The relative error assesses the other parameters:
\begin{equation}\label{eqn: model value error}
\frac{\left|\mathrm{value}_{\mathrm{sim}} - \mathrm{value}_{\mathrm{model}}\right|}{\mathrm{value}_{\mathrm{sim}}} \times 100.
\end{equation}

To construct stress cases for Assumption 9, which assumes an exponential distribution for the OOS lead time, we performed simulations using alternative distributions for lead time with the same MTTR, $\rrate{\oos}^{-1}$, but different coefficients of variation (CVs). Here, CV is defined as the ratio of the standard deviation to the mean. The tested cases were deterministic, gamma, and lognormal. For the gamma distribution, we considered CV values of $1/2$, $1/\sqrt{2}$, $1$, $\sqrt{2}$, and $2$, whereas for the lognormal distribution, we considered CV values of $\sqrt{2}$, $2$, and $4$. The deterministic case corresponds to a fixed lead time equal to $\rrate{\oos}^{-1}$. The gamma distribution was included because it generalizes the exponential distribution, with CV $=1$ corresponding to the baseline exponential case (Assumption 9). Therefore, both lower- and higher-variability cases relative to this baseline were tested. The lognormal distribution was included to examine long-tail stress cases, with CV values exceeding the baseline value of $1$. The results are summarized in Tables \ref{tab: average errors stress 98}--\ref{tab: average errors stress 90}, corresponding to the cases where $\ofr{\pln}^{\mathrm{req}}=\ofr{\prk}^{\mathrm{req}}=0.98$, $0.95$, and $0.90$, respectively.

\begin{table}[hbt!]
    \centering
    \rev{
    \caption{\rev{Averaged errors of result parameters, $\ofr{\pln}^{\mathrm{req}}=\ofr{\prk}^{\mathrm{req}}=0.98$}}
    \label{tab: average errors stress 98}
    \begin{tabular}{lcccccccccc}
        \hline\hline
        & Deterministic & \multicolumn{5}{c}{Gamma} & \multicolumn{3}{c}{Lognormal} \\
        \cmidrule(lr){2-2} \cmidrule(lr){3-7} \cmidrule(lr){8-10}
        $\mathrm{CV}$ & 0 & $1/2$ & $1/\sqrt{2}$ & $1$ (Exponential) & $\sqrt{2}$ & $2$ & $\sqrt{2}$ & $2$ & $4$ \\\hline
        \textbf{Relative Errors (\%)} & & & & & & & & & \\
        $\stklv\pln$ & 1.89 & 1.93 & 1.92 & 1.90 & 1.96 & 1.91 & 1.89 & 1.86 & 1.87 \\
        $\stklv\prk$ & 0.73 & 0.76 & 0.75 & 0.66 & 0.76 & 0.70 & 0.75 & 0.76 & 0.82 \\
        $\stklv{\wait}$ & 1.08 & 1.95 & 2.01 & 2.06 & 2.51 & 2.75 & 2.01 & 2.29 & 3.04 \\
        $\ofq{\pln}$ & 0.54 & 0.50 & 0.52 & 0.56 & 0.45 & 0.47 & 0.55 & 0.54 & 0.50 \\
        $\ofq{\prk}$ & 1.51 & 1.71 & 1.64 & 1.40 & 1.46 & 1.38 & 1.54 & 1.72 & 1.53 \\
        $\ofq{\oos}$ & 0.64 & 0.76 & 0.75 & 0.70 & 0.76 & 0.70 & 0.66 & 0.77 & 0.67 \\
        $\mtm{d}$ & 0.43 & 0.42 & 0.46 & 0.39 & 0.44 & 0.46 & 0.48 & 0.48 & 0.48 \\
        $\AMC$ & 0.62 & 0.60 & 0.60 & 0.60 & 0.66 & 0.59 & 0.56 & 0.66 & 0.70 \\
        \addlinespace
        \textbf{Absolute Errors (\%p)} & & & & & & & & & \\
        $\ofr\pln$ & 0.10 & 0.10 & 0.10 & 0.09 & 0.10 & 0.11 & 0.10 & 0.11 & 0.10 \\
        $\ofr\prk$ & 0.07 & 0.06 & 0.05 & 0.05 & 0.06 & 0.07 & 0.07 & 0.06 & 0.06 \\
        \hline\hline
    \end{tabular}
    }
\end{table}

\begin{table}[hbt!]
    \centering
    \rev{
    \caption{\rev{Averaged errors of result parameters, $\ofr{\pln}^{\mathrm{req}}=\ofr{\prk}^{\mathrm{req}}=0.95$}}
    \label{tab: average errors stress 95}
    \begin{tabular}{lcccccccccc}
        \hline\hline
        & Deterministic & \multicolumn{5}{c}{Gamma} & \multicolumn{3}{c}{Lognormal} \\
        \cmidrule(lr){2-2} \cmidrule(lr){3-7} \cmidrule(lr){8-10}
        $\mathrm{CV}$ & 0 & $1/2$ & $1/\sqrt{2}$ & $1$ (Exponential) & $\sqrt{2}$ & $2$ & $\sqrt{2}$ & $2$ & $4$ \\\hline
        \textbf{Relative Errors (\%)} & & & & & & & & & \\
        $\stklv\pln$ & 2.65 & 2.71 & 2.71 & 2.65 & 2.72 & 2.72 & 2.69 & 2.64 & 2.67 \\
        $\stklv\prk$ & 0.84 & 0.82 & 0.85 & 0.81 & 0.88 & 0.76 & 0.87 & 0.80 & 0.89 \\
        $\stklv{\wait}$ & 1.18 & 1.83 & 1.92 & 2.03 & 2.35 & 2.53 & 1.98 & 2.21 & 2.95 \\
        $\ofq{\pln}$ & 0.55 & 0.51 & 0.49 & 0.54 & 0.44 & 0.47 & 0.56 & 0.51 & 0.47 \\
        $\ofq{\prk}$ & 1.38 & 1.51 & 1.52 & 1.33 & 1.39 & 1.37 & 1.38 & 1.63 & 1.38 \\
        $\ofq{\oos}$ & 0.69 & 0.78 & 0.75 & 0.70 & 0.77 & 0.75 & 0.66 & 0.77 & 0.66 \\
        $\mtm{d}$ & 0.44 & 0.44 & 0.46 & 0.42 & 0.46 & 0.44 & 0.47 & 0.42 & 0.45 \\
        $\AMC$ & 0.71 & 0.64 & 0.69 & 0.63 & 0.72 & 0.63 & 0.64 & 0.69 & 0.71 \\
        \addlinespace
        \textbf{Absolute Errors (\%p)} & & & & & & & & & \\
        $\ofr\pln$ & 0.34 & 0.35 & 0.36 & 0.33 & 0.35 & 0.35 & 0.35 & 0.34 & 0.32 \\
        $\ofr\prk$ & 0.10 & 0.10 & 0.09 & 0.10 & 0.09 & 0.08 & 0.10 & 0.09 & 0.11 \\
        \hline\hline
    \end{tabular}
    }
\end{table}

\begin{table}[hbt!]
    \centering
    \rev{
    \caption{\rev{Averaged errors of result parameters, $\ofr{\pln}^{\mathrm{req}}=\ofr{\prk}^{\mathrm{req}}=0.90$}}
    \label{tab: average errors stress 90}
    \begin{tabular}{lcccccccccc}
        \hline\hline
        & Deterministic & \multicolumn{5}{c}{Gamma} & \multicolumn{3}{c}{Lognormal} \\
        \cmidrule(lr){2-2} \cmidrule(lr){3-7} \cmidrule(lr){8-10}
        $\mathrm{CV}$ & 0 & $1/2$ & $1/\sqrt{2}$ & $1$ (Exponential) & $\sqrt{2}$ & $2$ & $\sqrt{2}$ & $2$ & $4$ \\\hline
        \textbf{Relative Errors (\%)} & & & & & & & & & \\
        $\stklv\pln$ & 2.74 & 2.83 & 2.83 & 2.73 & 2.78 & 2.81 & 2.80 & 2.72 & 2.74 \\
        $\stklv\prk$ & 0.92 & 0.83 & 0.93 & 0.85 & 0.97 & 0.86 & 0.87 & 0.84 & 0.93 \\
        $\stklv{\wait}$ & 1.17 & 1.95 & 2.13 & 2.21 & 2.57 & 2.74 & 2.08 & 2.30 & 3.19 \\
        $\ofq{\pln}$ & 0.55 & 0.51 & 0.56 & 0.53 & 0.49 & 0.52 & 0.59 & 0.56 & 0.53 \\
        $\ofq{\prk}$ & 1.42 & 1.44 & 1.57 & 1.32 & 1.40 & 1.27 & 1.25 & 1.61 & 1.29 \\
        $\ofq{\oos}$ & 0.73 & 0.82 & 0.79 & 0.70 & 0.82 & 0.76 & 0.71 & 0.80 & 0.72 \\
        $\mtm{d}$ & 0.45 & 0.46 & 0.48 & 0.43 & 0.49 & 0.49 & 0.48 & 0.45 & 0.48 \\
        $\AMC$ & 0.74 & 0.67 & 0.70 & 0.65 & 0.71 & 0.66 & 0.67 & 0.71 & 0.72 \\
        \addlinespace
        \textbf{Absolute Errors (\%p)} & & & & & & & & & \\
        $\ofr\pln$ & 0.43 & 0.44 & 0.44 & 0.40 & 0.42 & 0.42 & 0.45 & 0.41 & 0.39 \\
        $\ofr\prk$ & 0.16 & 0.13 & 0.14 & 0.12 & 0.12 & 0.13 & 0.15 & 0.13 & 0.13 \\
        \hline\hline
    \end{tabular}
    }
\end{table}

The results show that the proposed model remains a good approximation across all tested cases. In the baseline case corresponding to Assumption 9, i.e., the gamma distribution with CV $=1$ (equivalently, the exponential distribution), all reported relative errors remain below $3$\% for the three service level requirements, and both absolute order fill rate errors remain below $0.5$\%p. This indicates that the model reproduces the simulation results well under the nominal assumption.

Among the reported parameters, $\stklv{\pln}$ and $\stklv{\wait}$ exhibit the largest errors. This is expected because these quantities are most directly coupled to the OOS-enabled closed-loop recovery process. However, their sensitivities to lead-time variability differ. In particular, $\stklv{\wait}$ is the most sensitive to the shape of the lead-time distribution: for both the gamma and lognormal families, its error increases as the CV increases, even though the MTTR is held fixed. This behavior is consistent with the fact that $\stklv{\wait}$ depends most directly on the residence time in the OOS stage. By contrast, the error in $\stklv{\pln}$ is relatively insensitive to the CV because the in-plane spare stock is influenced not only by recovery through OOS but also by replenishment from parking spares.

As the required order fill rate is relaxed from $0.98$ to $0.90$, the overall error level increases modestly, but the deterioration remains limited. Even in the most severe stress case---a lognormal lead time distribution with CV $=4$ under $\ofr{\pln}^{\mathrm{req}}=\ofr{\prk}^{\mathrm{req}}=0.90$---the largest average relative error is only $3.19$\%, observed for $\stklv{\wait}$, while the largest average absolute error in order fill rate across all tested cases is only $0.45$\%p. Therefore, although higher-order features of the lead time distribution do affect the prediction of the waiting stock, the resulting loss of accuracy remains small for the tested service-level requirements, including thresholds as low as $0.90$.

In summary, these results support the use of Assumption 9 as a tractable approximation for system-level modeling. Model accuracy is driven primarily by the mean OOS lead time, whereas the exact distributional form plays a secondary role. Accordingly, the exponential assumption provides a reasonable balance between analytical simplicity and predictive accuracy for design-space exploration and managerial trade-off analysis, although some caution is warranted when interpreting $\stklv{\wait}$ under highly variable, long-tailed lead-time distributions and low service-level requirements.}

\section{Decision-Making Scenarios and Optimization Problem Formulations}\label{section 4}
\rev{This study considers two primary stakeholders: the constellation operator and the OOS provider. The constellation operator aims to maintain the required service level of the constellation while minimizing maintenance cost. The OOS provider, in turn, aims to maximize profits by selecting service price and performance levels---in this study, responsiveness---that remain commercially attractive to the constellation operator. The decisions of the two stakeholders are strongly coupled: the operator's maintenance strategy determines the demand for OOS, whereas the provider's pricing and responsiveness influence the operator's preferred maintenance strategy. Balancing these tensions and designing appropriate incentives to achieve mutually beneficial outcomes that support the broader space industry are important challenges. Indeed, the interdependence between satellite operators and OOS providers is sometimes regarded as a chicken-and-egg problem \cite{us2025space}.

We first represent this interaction through the following bi-objective optimization problem, whose components are detailed in this section:
\newline\noindent
($\mathbf{P}$) Bi-Objective Optimization for the Constellation Operator and OOS Provider
\begin{equation}\label{eqn: P obj1}
    \min\ \AMC\left(\dec{};\prm{}\right),
\end{equation}
\begin{equation}\label{eqn: P obj2}
    \max\ \AP\left(\dec{};\prm{}\right),
\end{equation}
\noindent subject to
\begin{center}
    Eqs. \eqref{eqn: val cond 1}--\eqref{eqn: val cond 5},
\end{center}
\begin{equation}\label{eqn: P const 1}
    Qk_Q\le Q_{\max},
\end{equation}
\begin{equation}\label{eqn: P const 2}
    \scost{\oos}\left(\rrate{\oos};\prm{}\right)\le p_{\oos},
\end{equation}
\begin{equation}\label{eqn: P const 3}
    \AMC\left(\dec{};\prm{}\right)\le \AMC_{\mathrm{ref}},
\end{equation}
\begin{equation}\label{eqn: P const 4}
    \dec{}=\left[s,Q,k_s,k_Q,\nbr{\oos},\nbr{\prk},\alt{\prk},p_{\oos},\rrate{\oos}\right]
    \in\mathbb{Z}_{++}^{6}\times\mathcal{H}_{\prk}\times\mathbb{R}_{++}^{2}.
\end{equation}
Here, $\dec{}$ denotes the unified decision vector, and $\prm{}$ collects all remaining parameters. The decision vector contains both the constellation operator's decision variables, $\left[s,Q,k_s,k_Q,\nbr{\oos},\nbr{\prk},\alt{\prk}\right]$, and the OOS provider's decision variables, $\left[p_{\oos},\rrate{\oos}\right]$. In addition, $\mathcal{H}_{\prk}$ denotes the set of available parking orbit altitudes.

From the constellation operator's perspective, Eq. \eqref{eqn: P obj1} minimizes the annual maintenance cost. The service-level requirements in Eqs. \eqref{eqn: val cond 1} and \eqref{eqn: val cond 2} ensure that the maintenance system achieves the required order fill rates for in-plane and parking spares, respectively. Equation \eqref{eqn: P const 1} limits the number of loaded satellites so that it does not exceed the launch capacity, $Q_{\max}$, while Eqs. \eqref{eqn: val cond 3}--\eqref{eqn: val cond 5} impose the validity conditions of the analytical model.

From the OOS provider's perspective, Eq. \eqref{eqn: P obj2} represents the maximization of annual profit, defined as
\begin{equation}\label{eqn: oos profit}
    \AP=\left(p_{\oos}-\scost{\oos}\right)\ofq{\oos}\nbr{\pln},
\end{equation}
where $p_{\oos}$ and $\scost{\oos}$ are the unit price and cost of an OOS, respectively. The term $\ofq{\oos}\nbr{\pln}$ represents the total annual number of services performed, computed as the product of the OOS frequency, $\ofq{\oos}$, and the number of target orbits, $\nbr{\pln}$. Thus, annual profit equals the profit per service multiplied by the total annual service volume.

Given the OOS technologies, service types, and target subsystems, higher responsiveness generally requires greater operational effort, such as additional servicers, depots, and propellant for faster orbital transfers. To capture this trade-off, we adopt the following first-cut model relating cost and responsiveness:
\begin{equation}\label{eqn: cost responsiveness function}
\scost{\oos} = \scost{\oos,\min} + \frac{\alpha_1}{\left(\rrate{\oos}^{-1} - \rrate{\oos,\mathrm{ideal}}^{-1}\right)^{\alpha_2}}.
\end{equation}
Here, $\scost{\oos,\min}$ ($>0$) denotes the lower bound of the service cost, $\rrate{\oos,\mathrm{ideal}}$ ($>0$) denotes the upper bound of service responsiveness, and $\alpha_1$ ($>0$) and $\alpha_2$ ($>0$) are shape parameters. This function is strictly decreasing with respect to MTTR and exhibits diminishing returns as responsiveness approaches its ideal bound. Such a relationship, in which greater operational effort yields improved responsiveness at a diminishing rate, is consistent with empirical observations in terrestrial logistics \cite{chopra2019supply} and prior OOS-related studies \cite{dujonchay2017quantification, luu2020valuation, ho2020semi}.

\begin{figure}[hbt!]
    \centering
    \includegraphics[page=4, width=1.0\textwidth]{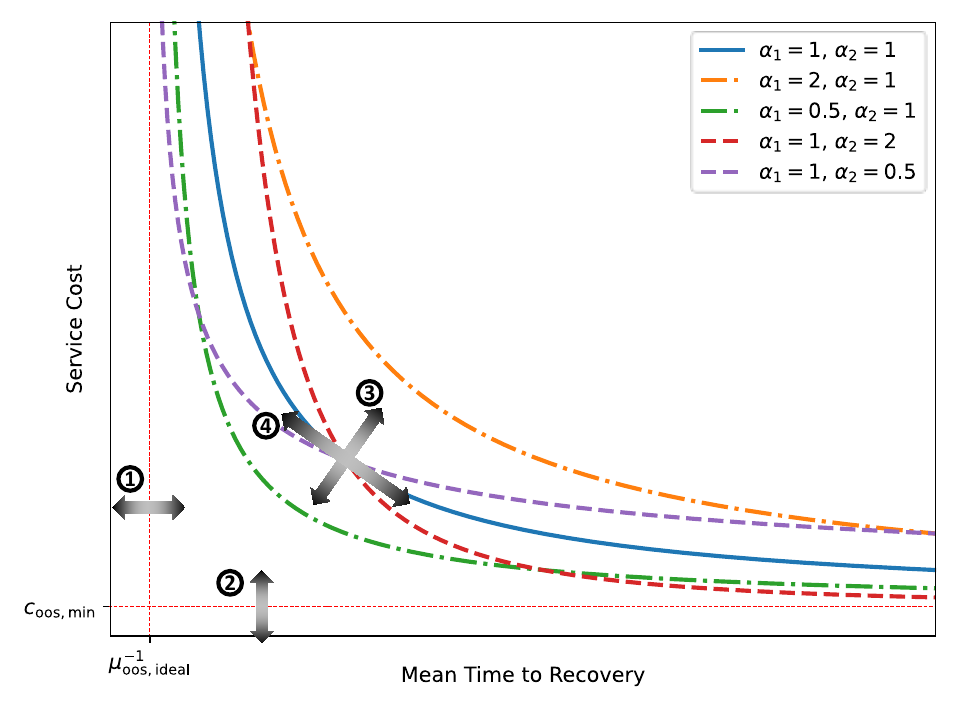}
    \caption{Visualization of the adopted cost-responsiveness function with the effect of each parameter (1: $\scost{\oos,\min}$, 2: $\rrate{\oos,\mathrm{ideal}}$, 3: $\alpha_1$, 4: $\alpha_2$).}
    \label{fig: oos cost structure}
\end{figure}

Equation \eqref{eqn: P const 2} ensures that the OOS price is no less than the corresponding service cost. Equation \eqref{eqn: P const 3} constrains the provider to choose the price and responsiveness so that the resulting OOS-integrated maintenance strategy remains acceptable to the constellation operator. The reference value $\AMC_{\mathrm{ref}}$ may be interpreted as the maintenance cost of a benchmark strategy without OOS, or as a premium-adjusted reference value reflecting the additional complexity and uncertainty introduced by OOS integration.

Solving Problem $\mathbf{P}$ yields a family of non-dominated solutions\cite{marler2004survey}, known as the Pareto front \cite{edgeworth1881mathematical, pareto2014manual}. In this study, the Pareto front is obtained using the non-dominated sorting genetic algorithm II (NSGA-II), a genetic algorithm-based metaheuristic for multi-objective optimization \cite{deb2002fast}.
}

\section{Case Study}\label{section 5}
This section presents a case study of the proposed bi-objective optimization problem \rev{($\mathbf{P}$)}. The case study demonstrates the applicability of the proposed framework using a realistically sized constellation and representative satellite specifications. We first establish and analyze a baseline scenario, then extend the analysis through a parametric study by varying the OOS-related parameters.

\subsection{Baseline Scenario}
We established the baseline scenario using the constellation's system parameters, launch service parameters, and OOS parameters, which are summarized in Tables \ref{tab: baseline parameters 1}, \ref{tab: baseline parameters 2}, and \ref{tab: baseline parameters 3}, respectively. \rev{The decision variables and their corresponding ranges are listed in Table \ref{tab: baseline decision variables}. The constellation and launch service parameters in Tables \ref{tab: baseline parameters 1} and \ref{tab: baseline parameters 2} are adopted from prior studies \cite{jakob2019optimal,kim2025replenishment} to reflect realistic, real-world-scale systems. By contrast, the OOS-related parameters cannot yet be assigned fully realistic values, because the relevant service capabilities and market conditions remain uncertain. Instead, they are specified within plausible ranges based on the existing literature, and their effects on the results are further examined in the parametric study by varying their values. In particular, the minimum MTTR is set to $2$ weeks (i.e., $\rrate{\oos}=1/2$), which lies within a reasonable range according to prior studies on OOS responsiveness \cite{dujonchay2017quantification, luu2020valuation, ho2020semi}. The serviceable failure ratio, $r_\oos$, is set to $25$\%, considering the maximum fraction of failures that could be induced by a single spacecraft subsystem failure, namely the thruster subsystem, as discussed by Saleh and Castet \cite{saleh2011spacecraft}. The minimum OOS cost is set equal to the satellite production cost, $\scost{\sat}$. The shape parameters, $\alpha_1$ and $\alpha_2$, are both set to 1.}

\begin{table}[hbt!]
    \centering
    \caption{Constellation's system parameters for baseline scenario}
    \begin{tabular}{lcc}
        \hline\hline
        Parameter, Unit & Notation & Value \\\hline
        Inclination, deg & $\inc$ & 60 \\
        Satellite failure rate, failures/year & $\frate{\sat}$ & 0.2 \\
        Number of orbital planes, planes & $\nbr{\pln}$ & 40 \\
        Number of satellites per orbital plane, satellites/plane & $\nbr{\sat}$ & 40 \\
        Altitude of orbital planes, km & $\alt{\pln}$ & 1200 \\
        Fuel mass conversion coefficient, \$M/kg & $\epsilon_{\mathrm{fuel}}$ & 0.01 \\
        Annual satellite holding cost, \$M/satellite/year & $\scost{\hld}$ & 0.5 \\
        Specific impulse of satellite's propulsion system, sec & $I_{\mathrm{sp}}$ & 1200 \\
        Dry mass, kg & $M_{\mathrm{dry}}$ & 150 \\
        Propellant mass rate, kg/s & $\dot{M}_{\mathrm{fuel}}$ & 1.3E-6 \\
        Production cost, \$M/satellite & $\scost{\sat}$ & 0.5 \\
        Lifespace of spacecraft, years & $\mtm{\mathrm{life}}$ & 30 \\
        \multirow{2}{*}{Service level requirement, -} & $\ofr{\pln}^{\mathrm{req}}$ & 0.98 \\
         & $\ofr{\prk}^{\mathrm{req}}$ & 0.98 \\
         Annual maintenance cost threshold for adopting OOS, \$M/year & $\AMC_{\mathrm{ref}}$ & 925.1 \\
        \hline\hline
    \end{tabular}
    \label{tab: baseline parameters 1}
\end{table}

\begin{table}[hbt!]
    \centering
    \caption{Launch service parameters for baseline scenario}
    \begin{tabular}{lcc}
        \hline\hline
        Parameter, Unit & Notation & Value \\\hline
        Launch cost, \$M & $\scost{\lau}$ & 67 \\
        Launch capacity, satellites & $Q_{\max}$ & 40 \\
        Launch order processing time, weeks & $t_{\lau}$ & 12 \\
        Mean waiting time to launch, weeks & $\rrate{\lau}^{-1}$ & 8 \\
        \hline\hline
    \end{tabular}
    \label{tab: baseline parameters 2}
\end{table}

\begin{table}[hbt!]
    \centering
    \caption{OOS parameters for baseline scenario}
    \begin{tabular}{lcc}
        \hline\hline
        Parameter, Unit & Notation & Value \\\hline
        Serviceable failures ratio, - & $r_{\oos}$ & 0.25 \\
        The upper bound of responsiveness parameter, week$^{-1}$ & $\rrate{\oos,\mathrm{ideal}}$ & 0.5 \\
        The lower bound of service cost, \$M & $\scost{\oos,\min}$ & 0.5 \\
        \multirow{2}{*}{Shape parameters for cost-reponsiveness function, -} & $\alpha_1$ & 1 \\
        & $\alpha_2$ & 1 \\
        \hline\hline
    \end{tabular}
    \label{tab: baseline parameters 3}
\end{table}

\begin{table}[hbt!]
    \centering
    \caption{Decision variables and their ranges for baseline scenario}
    \begin{tabular}{lcc}
        \hline\hline
        Decision Variable, Unit & Notation & Range \\\hline
        Reorder point of in-plane spares, satellites & $s$ & $\left[1, 20\right]$ \\
        Order quantity of in-plane spares, satellites & $Q$ & $\left[1, Q_{\max}\right]$ \\
        Reorder point of parking spares, - & $k_s$ & $\left[1, 20\right]$ \\
        Order quantity of parking spares, - & $k_Q$ & $\left[1, Q_{\max}\right]$ \\
        Maximum number of services, - & $\nbr{\oos}$ & $\left[1, 4\right]$ \\
        Number of parking orbits, planes & $\nbr{\prk}$ & $\left[1, 20\right]$ \\
        Altitude of parking orbits, km & $\alt{\prk}$ & $\left[500, 1000\right]$ \\
        OOS price, \$M & $p_{\oos}$ & $\left[\scost{\oos,\min},\scost{\oos,\min}+5\right]$ \\
        OOS MTTR, weeks & $\rrate{\oos}^{-1}$ & $\left(\rrate{\oos,{\mathrm{ideal}}}^{-1},12\right]$ \\
        \hline\hline
    \end{tabular}
    \label{tab: baseline decision variables}
\end{table}

\rev{In addition, the reference annual maintenance cost ($\AMC_{\mathrm{ref}}$) in Table \ref{tab: baseline parameters 1} is derived from the baseline multi-echelon spare strategy proposed by Jakob et al. \cite{jakob2019optimal}, which has the same supply chain structure as that considered in this study but does not integrate OOS. The optimization uses the same constellation and launch service parameters as those listed in Tables \ref{tab: baseline parameters 1} and \ref{tab: baseline parameters 2}, and the resulting optimal solution is presented in Table \ref{tab: baseline solutions without OOS}. Among the resulting outputs, only the optimal annual maintenance cost, \$925.1M, is used in the baseline scenario and adopted as the reference annual maintenance cost, $\AMC_{\mathrm{ref}}$, in Table \ref{tab: baseline parameters 1}. This solution was obtained using a genetic algorithm implemented in Python with the \texttt{pymoo} package \cite{blank2020pymoo}, with a population size of 400 and 200 generations, while all other algorithm settings were left at their default values. The optimization was performed over the same ranges of the decision variables given in Table \ref{tab: baseline decision variables}, except for the OOS-related variables---$p_\oos$ and $\rrate{\oos}$---which are not relevant to the multi-echelon spare strategy without OOS.}

\begin{table}[hbt!]
    \centering
    \caption{Optimal solution and result parameters of the baseline scenario without OOS}
    \begin{tabular}{lc}
        \hline\hline
        Result & Value \\\hline
        \textbf{Decision Variables} &\\
        $s$, satellites & 4 \\
        $Q$, satellites & 4 \\
        $k_s$, - & 10 \\
        $k_Q$, - & 10 \\
        $\nbr{\prk}$, planes & 6 \\
        $\alt{\prk}$, km & 795.4 \\\addlinespace
        \textbf{Performance Outputs} &\\
        $\AMC$, \$M/year & 925.1\\
        $\acost{\lau}$, \$M/year & 536.0\\
        $\acost{\hld}$, \$M/year & 220.8\\
        $\acost{\mnf}$, \$M/year & 160.0\\
        $\ofr{\pln}$, \% & 98.0\\
        $\ofr{\prk}$ \% & 98.3\\
        \hline\hline
    \end{tabular}
    \label{tab: baseline solutions without OOS}
\end{table}

\begin{figure}[hbt!]
    \centering
    \includegraphics[page=5, width=1.0\textwidth]{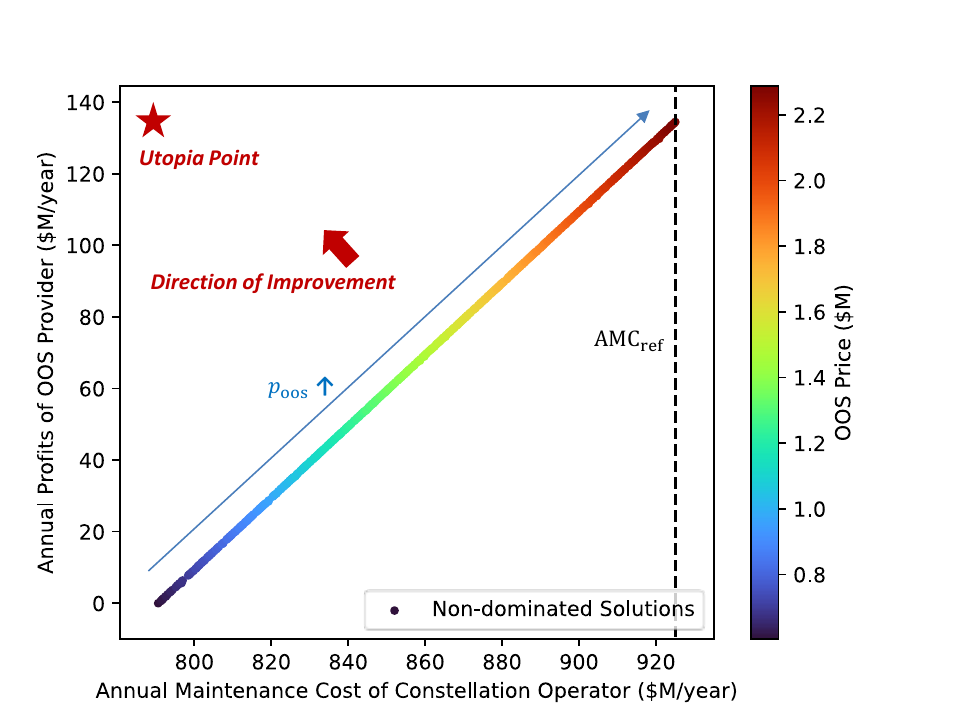}
    \caption{Solutions of the baseline scenario scattered in the objective space (color: OOS price).}
    \label{fig: baseline solutions}
\end{figure}

\rev{Using these parameter settings, the bi-objective optimization problem ($\mathbf{P}$) was solved using NSGA-II \cite{deb2002fast}, implemented in Python with the \texttt{pymoo} package \cite{blank2020pymoo}. The algorithm was run with a population size of 400 for 200 generations, while all other algorithm parameters were kept at their default values. To assess the reliability of the results, we performed three independent runs with different random seeds. Because the ranges of the obtained solutions were sufficiently similar across the three runs in both the decision and objective spaces, we retained the solution set generated using the first random seed.} The obtained optimal solutions are scattered in the objective space in Fig. \ref{fig: baseline solutions}. The $x$-axis and $y$-axis in Fig. \ref{fig: baseline solutions} represent the first objective (i.e., the constellation operator's annual maintenance cost, $\AMC$) and the second objective (i.e., the OOS provider's annual profit $\AP$), respectively. Since the first objective is minimized and the second is maximized, the utopia point is located at the upper-left corner of the figure, marked by a red star.

The values and ranges of the decision variables, together with the resulting performance measures of the Pareto-optimal solutions, are summarized in Table \ref{tab: baseline solutions}. \rev{Figure \ref{fig: baseline solutions} shows that, under the baseline setting, the Pareto-optimal solutions form an almost linear frontier spanning $\AMC=790.6$--$925.0$ \$/year and $\AP=0$--$134.5$ \$/year. This frontier is bounded by the OOS provider's non-negative profit condition (Eq. \eqref{eqn: P const 2}) and the constellation operator's annual maintenance cost threshold, $\AMC_{\mathrm{ref}}$ (Eq. \eqref{eqn: P const 3}). In addition, the constraints in Eqs. \eqref{eqn: val cond 1}--\eqref{eqn: val cond 3} and Eq. \eqref{eqn: P const 1} are effectively binding along the Pareto frontier. Although they are not satisfied at equality for every solution, the remaining slack is small and is attributable to the integrality of the decision variables. Relative to the benchmark solution without OOS in Table \ref{tab: baseline solutions without OOS}, the constellation operator can reduce annual maintenance cost by as much as \$134.5M/year (14.5\%) while maintaining the required service levels. Moreover, the Pareto-optimal set is generated by an almost invariant maintenance strategy on the constellation operator side: $s=3$, $Q=4$, $k_s=6$, $k_Q=10$, $\nbr{\prk}=7$, $\alt{\prk}\in[674.3,700.4]$ km, and $\nbr{\oos}=4$. This indicates that, in the baseline environment, the location of a Pareto solution is governed primarily by the OOS provider's decision rather than by multiple qualitatively different maintenance strategies for the constellation operator.

Compared with the optimal strategy without OOS, the OOS-integrated solutions reduce annual launch cost from \$536.0M/year to \$402.4M/year and annual manufacturing cost from \$160.0M/year to \$120.1M/year, both reductions of approximately 25\%. These savings are partially offset by annual OOS expenditures of \$48.1M/year to \$182.4M/year. Since $\gamma_0=75.1\%$, approximately 24.9\% of failures are handled through OOS rather than by newly launched replacements. Along the Pareto frontier, the main variation comes from the OOS price $p_{\oos}$, whereas the optimal OOS MTTR remains near its upper bound, $\rrate{\oos}^{-1}=10.6$--$12.0$ weeks. This implies that, under the assumed cost-responsiveness relationship, the provider maximizes profit by selecting relatively slow but inexpensive service. Accordingly, the Pareto frontier is nearly linear, indicating that movement along it primarily redistributes value between the constellation operator and the OOS provider, rather than substantially altering the physical structure of the maintenance system.}

\begin{table}[hbt!]
    \centering
    \caption{Decision variables and resultant parameters of the Pareto-optimal solutions: baseline scenario}
    \begin{tabular}{lc}
        \hline\hline
        Result, Unit & Value or Range \\\hline
        \textbf{Decision Variables} & \\
        $s$, satellites & 3 \\
        $Q$, satellites & 4 \\
        $k_s$, - & 6 \\
        $k_Q$, - & 10 \\
        $\nbr{\prk}$, planes & 7 \\
        $\alt{\prk}$, km & [674.3, 700.4] \\
        $\nbr{\oos}$, - & 4 \\
        $p_{\oos}$, \$M & [0.6, 2.3] \\
        $\rrate{\oos}^{-1}$, weeks & [10.6, 12.0] \\\addlinespace
        \textbf{Performance Outputs} & \\
        $\AMC$, \$M/year & [790.6, 925.0]\\
        $\AP$, \$M/year & [0, 134.5]\\
        $\acost{\lau}$, \$M/year & 402.4\\
        $\acost{\hld}$, \$M/year & [211.4, 212.6]\\
        $\acost{\mnf}$, \$M/year & 120.1\\
        $\acost{\oos}$, \$M/year & [48.1, 182.4]\\
        $\mtm{d}/\nbr{t}$, years & [29.2, 29.3]\\
        $\ofr{\pln}$, \% & [98.0, 98.1]\\
        $\ofr{\prk}$, \% & 98.1\\
        $\gamma_0$, \% & 75.1\\
        \hline\hline
    \end{tabular}
    \label{tab: baseline solutions}
\end{table}

\subsection{Parametric Study}
\rev{We next examine how the Pareto-optimal solutions change when the OOS-related parameters are perturbed around the baseline scenario. Specifically, we vary the serviceable-failure ratio $r_{\oos}$, the lower bound of OOS cost $\scost{\oos,\min}$, the ideal responsiveness parameter $\rrate{\oos,\mathrm{ideal}}$, and the two shape parameters $\alpha_1$ and $\alpha_2$, as summarized in Table \ref{tab: parametric parameters}. All optimization runs were performed using the same algorithm settings as in the baseline scenario. To avoid repetitive interpretation, the results are discussed in terms of three questions: how each parameter shifts the Pareto frontier, whether it changes the constellation operator's preferred maintenance strategy, and what mechanism explains the observed change.}

\begin{table}[hbt!]
    \centering
    \caption{Cost-responsiveness function parameters for parametric study instances}
    \begin{tabular}{lccccc}
        \hline\hline
        & $r_{\oos}$ & $\scost{\oos,\min}$ & $\rrate{\oos,\mathrm{ideal}}$ & $\alpha_1$ & $\alpha_2$\\
        Instance ID & (-) & (\$M) & (week$^{-1}$) & (-) & (-)\\\hline
        Instance 0 (Baseline) & 0.25 & 0.5 & 0.5 & 1 & 1 \\
        Instance 1-1 & \textbf{0.5} & 0.5 & 0.5 & 1 & 1 \\
        Instance 1-1 & \textbf{0.1} & 0.5 & 0.5 & 1 & 1 \\
        Instance 2-1 & 0.25 & \textbf{1} & 0.5 & 1 & 1 \\
        Instance 2-2 & 0.25 & \textbf{0.25} & 0.5 & 1 & 1 \\
        Instance 3-1 & 0.25 & 0.5 & \textbf{1} & 1 & 1 \\
        Instance 3-2 & 0.25 & 0.5 & \textbf{0.25} & 1 & 1 \\
        Instance 4-1 & 0.25 & 0.5 & 0.5 & \textbf{2} & 1 \\
        Instance 4-2 & 0.25 & 0.5 & 0.5 & \textbf{0.5} & 1 \\
        Instance 5-1 & 0.25 & 0.5 & 0.5 & 1 & \textbf{2} \\
        Instance 5-2 & 0.25 & 0.5 & 0.5 & 1 & \textbf{0.5} \\
        \hline\hline
    \end{tabular}
    \label{tab: parametric parameters}
\end{table}

\subsubsection{Serviceable Failures Ratio: \texorpdfstring{$r_{\oos}$}{}}
\rev{Among the tested parameters, the serviceable-failure ratio $r_{\oos}$ has the strongest structural effect, because it directly determines the share of failures that can be covered by the production-based replenishment to the recovery from OOS. As shown in Fig. \ref{fig: parametric solutions 1} and Table \ref{tab: parametric solutions 1}, increasing $r_{\oos}$ from $0.1$ to $0.5$ shifts the Pareto frontier markedly toward the utopia point: the minimum $\AMC$ decreases from \$873.4M/year to \$661.9M/year, while the maximum $\AP$ increases from \$51.7M/year to \$263.1M/year. This improvement is accompanied by substantial reductions in annual launch and manufacturing costs, from \$482.4M/year to \$283.7M/year and from \$144.0M/year to \$82.6M/year, respectively. At the same time, $\gamma_0$ decreases from 90.0\% to 51.6\%, implying that the fraction of failures handled through OOS increases from 10.0\% to 48.4\%.

\begin{figure}[htpb]
    \centering
    \begin{subfigure}[b]{0.48\textwidth}
        \centering
        \includegraphics[page=6, width=0.99\textwidth]{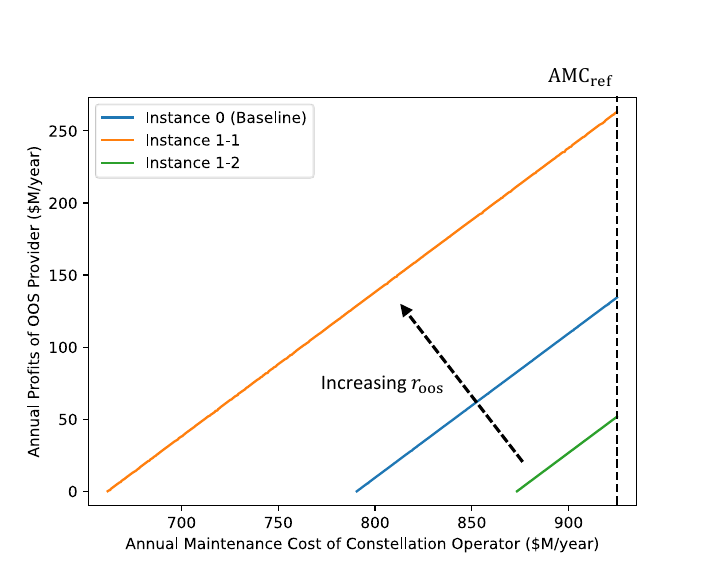}
        \caption{Instance 1-X: variations in $r_\oos$.}
        \label{fig: parametric solutions 1}
    \end{subfigure}
    \begin{subfigure}[b]{0.48\textwidth}
        \centering
        \includegraphics[page=7, width=0.99\textwidth]{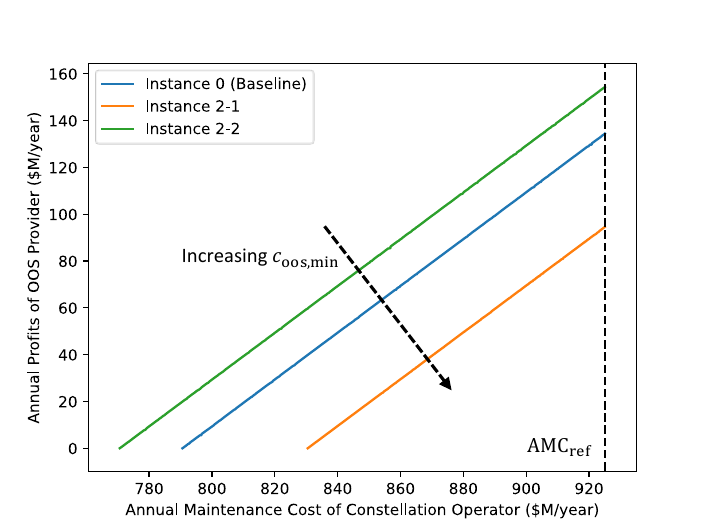}
        \caption{Instance 2-X: variations in $\scost{\oos,\min}$.}
        \label{fig: parametric solutions 2}
    \end{subfigure}
    \begin{subfigure}[b]{0.48\textwidth}
        \centering
        \includegraphics[page=8, width=0.99\textwidth]{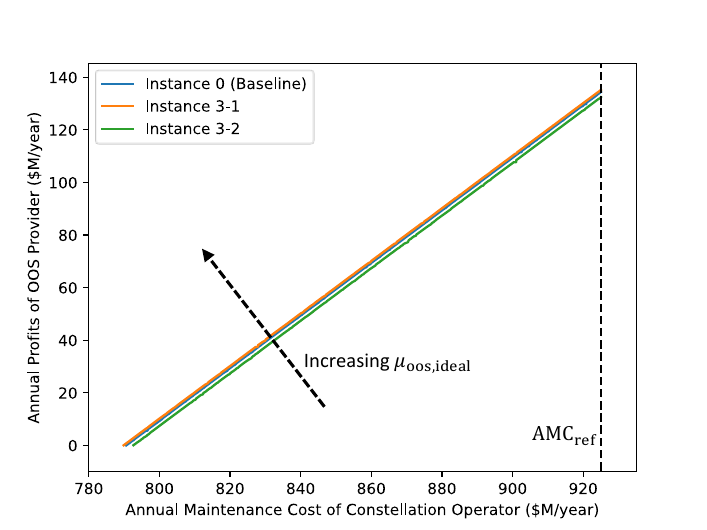}
        \caption{Instance 3-X: variations in $\rrate{\oos,\mathrm{ideal}}$.}
        \label{fig: parametric solutions 3}
    \end{subfigure}
    \begin{subfigure}[b]{0.48\textwidth}
        \centering
        \includegraphics[page=9, width=0.99\textwidth]{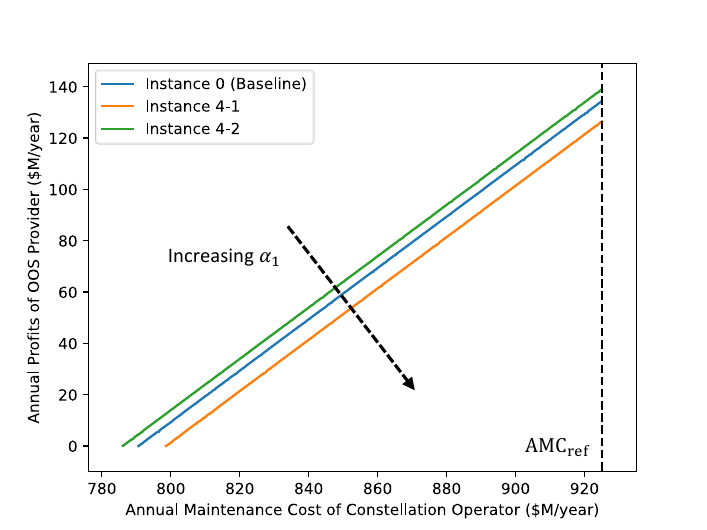}
        \caption{Instance 4-X: variations in $\alpha_1$.}
        \label{fig: parametric solutions 4}
    \end{subfigure}
    \begin{subfigure}[b]{0.48\textwidth}
        \centering
        \includegraphics[page=10, width=0.99\textwidth]{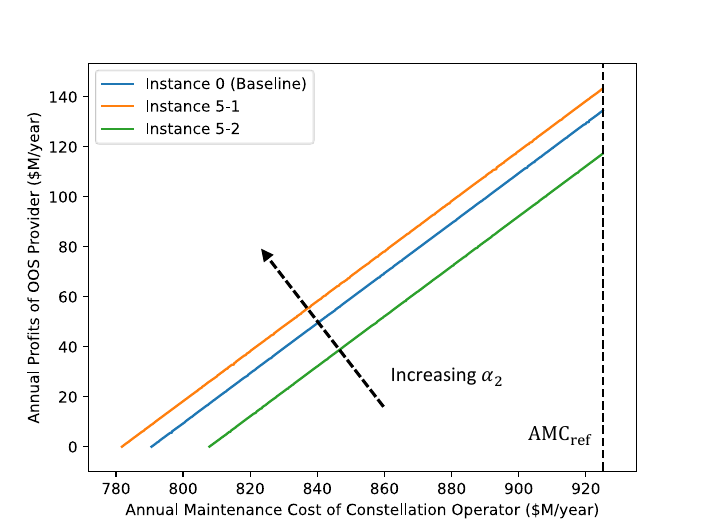}
        \caption{Instance 5-X: variations in $\alpha_2$.}
        \label{fig: parametric solutions 5}
    \end{subfigure}
    \caption{Pareto fronts of each Instance in the objective space.}
    \label{fig: parametric solutions}
\end{figure}

\begin{table}[hbt!]
    \centering
    \caption{Parametric study results: variations in $r_\oos$}
    \begin{tabular}{lcc}
        \hline\hline
        Result, Unit & Instance 1-1 & Instance 1-2\\
         & ($r_{\oos}=0.5$) & ($r_{\oos}=0.1$)\\\hline
        \textbf{Decision Variables} & & \\
        $s$, satellites & 3 & 4 \\
        $Q$, satellites & 3 & 4\\
        $k_s$, - & 5 & 9 \\
        $k_Q$, - & 13 & 10 \\
        $\nbr{\prk}$, planes & 7 & 6 \\
        $\alt{\prk}$, km & [638.3, 700.4] & [818.3, 823.3] \\
        $\nbr{\oos}$, - & 4 & [3, 4] \\
        $p_{\oos}$, \$M & [0.6, 2.3] & [0.6, 2.2] \\
        $\rrate{\oos}^{-1}$, weeks & [10.4, 12.0] & [10.4, 12.0]  \\\addlinespace
        \textbf{Performance Outputs} & & \\
        $\AMC$, \$M/year & [661.9, 925.0] & [873.4, 925.0] \\
        $\AP$, \$M/year & [0, 263.1] & [0, 51.7] \\
        $\acost{\lau}$, \$M/year & 283.7 & 482.4 \\
        $\acost{\hld}$, \$M/year & [194.9, 197.4] & [220.2, 220.9] \\
        $\acost{\mnf}$, \$M/year & 82.6 & 144.0 \\
        $\acost{\oos}$, \$M/year & [93.0, 356.3] & [19.4, 71.0] \\
        $\mtm{d}/\nbr{t}$, years & [29.2, 29.4] & [23.5, 29.4] \\
        $\ofr{\pln}$, \% & [98.0, 98.1] & 98.0 \\
        $\ofr{\prk}$, \% & 98.0 & 98.4 \\
        $\gamma_0$, \% & 51.6 & 90.0 \\
        \hline\hline
    \end{tabular}
    \label{tab: parametric solutions 1}
\end{table}

Unlike the other OOS-related parameters, $r_{\oos}$ also changes the constellation operator's preferred maintenance configuration. When $r_{\oos}=0.5$, the optimal value of $\nbr{\oos}$ remains fixed at its upper bound. In contrast, when $r_{\oos}=0.1$, $\nbr{\oos}$ takes values of 3 or 4, and the optimal stocking decisions move closer to those of the without-OOS benchmark in Table \ref{tab: baseline solutions without OOS}: $s$, $Q$, $k_Q$, and $\nbr{\prk}$ coincide with the benchmark values, while $k_s$ and $\alt{\prk}$ differ only modestly. This result indicates that $r_{\oos}$ is the key parameter determining whether OOS functions merely as an economic supplement or becomes a structurally important component of the maintenance architecture.}

\subsubsection{Shift Parameters of Cost-Responsiveness Function: \texorpdfstring{$\scost{\oos,\min}$ \& $\rrate{\oos,\mathrm{ideal}}$}{}}

\rev{We next examine the shift parameters of the cost-responsiveness function, $\scost{\oos,\min}$ and $\rrate{\oos,\mathrm{ideal}}$. The resulting solutions are illustrated in Figs. \ref{fig: parametric solutions 2} and \ref{fig: parametric solutions 3}, while the detailed values and ranges of the decision variables and outcome measures are reported in Table \ref{tab: parametric solutions 2 and 3}. Their main effect is to shift the Pareto frontier without structural change of the constellation operator's preferred maintenance strategy. Across Instances 2-X and 3-X, the operator-side decision variables remain identical to those of the baseline scenario, and the corresponding dependence on OOS also remains unchanged, as reflected by the constant values of $\acost{\lau}=\$402.4$M/year, $\acost{\mnf}=\$120.1$M/year, and $\gamma_0=75.1\%$.

\begin{table}[hbt!]
    \centering
    \caption{Parametric study results: variations in $\scost{\oos,\min}$ and $\rrate{\oos,\mathrm{ideal}}$}
    \begin{tabular}{lcccc}
        \hline\hline
        Result, Unit & Instance 2-1 & Instance 2-2 & Instance 3-1 & Instance 3-2\\
         & ($\scost{\oos,\min}=1$) & ($\scost{\oos,\min}=0.25$) & ($\rrate{\oos,\mathrm{ideal}}=1$) & ($\rrate{\oos,\mathrm{ideal}}=0.25$)\\\hline
        \textbf{Decision Variables} & & & & \\
        $s$, satellites & \multicolumn{4}{c}{3} \\
        $Q$, satellites & \multicolumn{4}{c}{4}\\
        $k_s$, - & \multicolumn{4}{c}{6}\\
        $k_Q$, - & \multicolumn{4}{c}{10}\\
        $\nbr{\prk}$, planes & \multicolumn{4}{c}{7}\\
        $\alt{\prk}$, km & [685.3, 700.4] & [675.5, 700.4] & [676.5, 700.4] & [656.5, 700.4]  \\
        $\nbr{\oos}$, - & \multicolumn{4}{c}{4} \\
        $p_{\oos}$, \$M & [1.1, 2.3] & [0.4, 2.3] & [0.6, 2.3] & [0.6, 2.3]  \\
        $\rrate{\oos}^{-1}$, weeks & [10.8, 12.0] & [10.3, 12.0] & [9.3, 12.0] & [11.2, 12.0] \\\addlinespace
        \textbf{Performance Outputs} & & & & \\
        $\AMC$, \$M/year & [830.5, 925.0] & [770.7, 925.0] & [789.9, 925.0] & [792.7, 925.0] \\
        $\AP$, \$M/year & [0, 94.6] & [0, 154.4] & [0.1, 135.2] & [0.1, 132.4] \\
        $\acost{\lau}$, \$M/year & \multicolumn{4}{c}{402.4}\\
        $\acost{\hld}$, \$M/year & [211.1, 212.6] & [211.1, 212.6] & [210.4, 212.4] & [211.8, 212.5] \\
        $\acost{\mnf}$, \$M/year & \multicolumn{4}{c}{120.1}\\
        $\acost{\oos}$, \$M/year & [87.8, 182.3] & [28.0, 182.5] & [47.8, 182.9] & [50.0, 182.3] \\
        $\mtm{d}/\nbr{t}$, years & [29.2, 29.3] & [29.1, 29.3] & [29.2, 29.3] & [29.2, 29.3] \\
        $\ofr{\pln}$, \% & 98.0 & [98.0, 98.1] & [98.0, 98.1] & [98.0, 98.1] \\
        $\ofr{\prk}$, \% & \multicolumn{4}{c}{98.1}\\
        $\gamma_0$, \% & \multicolumn{4}{c}{75.1}\\
        \hline\hline
    \end{tabular}
    \label{tab: parametric solutions 2 and 3}
\end{table}

The effect of $\scost{\oos,\min}$ is direct and substantial. Increasing $\scost{\oos,\min}$ from \$0.25M to \$1M shifts the Pareto frontier downward and to the right, increasing the minimum $\AMC$ from \$770.7M/year to \$830.5M/year and reducing the maximum $\AP$ from \$154.4M/year to \$94.6M/year. By contrast, varying $\rrate{\oos,\mathrm{ideal}}$ has only a minor effect: changing it from $0.25$ to $1$ week$^{-1}$ changes the minimum $\AMC$ only from \$792.7M/year to \$789.9M/year and the maximum $\AP$ only from \$132.4M/year to \$135.2M/year. This limited sensitivity arises because the efficient solutions remain in the low-responsiveness region of the cost-responsiveness curve, where horizontal shifts of the asymptote have only a small influence on cost.}

\subsubsection{Shape Parameters of Cost-Responsiveness Function: \texorpdfstring{$\alpha_1$ \& $\alpha_2$}{}}
\rev{Lastly, we vary the shape parameters $\alpha_1$ and $\alpha_2$, which reshape the curvature of the cost-responsiveness function rather than shift it rigidly. The resulting solutions are presented in Figs. \ref{fig: parametric solutions 4} \& \ref{fig: parametric solutions 5} and Table \ref{tab: parametric solutions 4 and 5}. Their effects are moderate and, again, do not alter the constellation operator's maintenance strategy. Across all four instances, the operator-side decision variables remain identical to those in the baseline scenario, and the values of $\acost{\lau}$, $\acost{\mnf}$, and $\gamma_0$ therefore remain unchanged.

\begin{table}[hbt!]
    \centering
    \caption{Parametric study results: variations in $\alpha_1$ and $\alpha_2$}
    \begin{tabular}{lcccc}
        \hline\hline
        Result, Unit & Instance 4-1 & Instance 4-2 & Instance 5-1 & Instance 5-2\\
         & ($\alpha_1=2$) & ($\alpha_1=0.5$) & ($\alpha_2=2$) & ($\alpha_2=0.5$)\\\hline
         \textbf{Decision Variables} & & & & \\
        $s$, satellites & \multicolumn{4}{c}{3} \\
        $Q$, satellites & \multicolumn{4}{c}{4}\\
        $k_s$, - & \multicolumn{4}{c}{6}\\
        $k_Q$, - & \multicolumn{4}{c}{10}\\
        $\nbr{\prk}$, planes & \multicolumn{4}{c}{7}\\
        $\alt{\prk}$, km & [675.1, 700.4] & [683.7, 700.4] & [672.1, 700.4] & [684.5, 700.4]  \\
        $\nbr{\oos}$, - & \multicolumn{4}{c}{4} \\
        $p_{\oos}$, \$M & [0.7, 2.3] & [0.6, 2.3] & [0.5, 2.3] & [0.8, 2.3]  \\
        $\rrate{\oos}^{-1}$, weeks & [11.7, 12.0] & [7.7, 10.9] & [6.8, 9.3] & [11.5, 12.0] \\\addlinespace
        \textbf{Result Parameters} & \\
        $\AMC$, \$M/year & [798.6, 925.0] & [786.1, 925.0] & [781.8, 925.0] & [807.9, 925.0] \\
        $\AP$, \$M/year & [0.0, 126.5] & [0.1, 139.0] & [0.0, 143.2] & [0.0, 117.2] \\
        $\acost{\lau}$, \$M/year & \multicolumn{4}{c}{402.4}\\
        $\acost{\hld}$, \$M/year & [212.2, 212.6] & [209.1, 211.6] & [208.4, 210.4] & [212.0, 212.6] \\
        $\acost{\mnf}$, \$M/year & \multicolumn{4}{c}{120.1}\\
        $\acost{\oos}$, \$M/year & [55.9, 182.3] & [45.1, 184.4] & [42.0, 185.5] & [65.2, 182.3] \\
        $\mtm{d}/\nbr{t}$, years & [29.2, 29.3] & [28.9, 29.2] & [28.9, 29.1] & [29.2, 29.3] \\
        $\ofr{\pln}$, \% & [98.0, 98.1] & 98.0 & [98.0, 98.1] & [98.0, 98.1] \\
        $\ofr{\prk}$, \% & \multicolumn{4}{c}{98.1}\\
        $\gamma_0$, \% & \multicolumn{4}{c}{75.1}\\
        \hline\hline
    \end{tabular}
    \label{tab: parametric solutions 4 and 5}
\end{table}

Increasing $\alpha_1$ makes the cost-responsiveness curve less favorable and shifts the Pareto frontier away from the utopia point: the minimum $\AMC$ increases from \$786.1M/year for $\alpha_1=0.5$ to \$798.6M/year for $\alpha_1=2$, while the maximum $\AP$ decreases from \$139.0M/year to \$126.5M/year. The effect of $\alpha_2$ is reversed in the present problem because the efficient solutions lie near the low-responsiveness tail of the curve. Specifically, increasing $\alpha_2$ from $0.5$ to $2$ improves the frontier, reducing the minimum $\AMC$ from \$807.9M/year to \$781.8M/year and increasing the maximum $\AP$ from \$117.2M/year to \$143.2M/year. Hence, the curvature of the provider's cost technology matters, but mainly through the economic terms of OOS rather than through a redesign of the constellation maintenance system.}

\subsection{Discussions}
\rev{
\subsubsection{Main Findings and Strategic Implications}
Three main findings emerge from the case study. First, under fixed constellation, launch service, and OOS-related parameters, the Pareto-optimal set within a given scenario (with fixed OOS-related parameters) is generated by an almost invariant maintenance strategy on the constellation operator's side. In the baseline scenario, the optimal operator-side decisions remain unchanged across the Pareto frontier, while the main variation comes from the provider-side pricing decision, with only limited variation in OOS responsiveness. As a result, the baseline Pareto frontier spans approximately $(\AMC,\AP)=(\$925.0\text{M/year},0)$ to $(\$790.6\text{M/year},\$134.5\text{M/year})$. This indicates that, relative to the benchmark solution without OOS, OOS can create up to \$134.5M/year of maintenance-cost savings for the constellation operator while still allowing non-negative profit for the OOS provider. A similar within-instance pattern is also observed across every scenario in the parametric study.

Second, OOS provides material logistical relief even before considering how value is divided between the two stakeholders. Relative to the benchmark solution without OOS, the baseline OOS-integrated solutions reduce annual launch cost from \$536.0M/year to \$402.4M/year and annual manufacturing cost from \$160.0M/year to \$120.1M/year, corresponding to reductions of approximately 25\% in both components, while maintaining the required service levels of approximately 98\% throughout the case study. In the baseline case, $\gamma_0=75.1\%$, implying that approximately 24.9\% of failures are handled through OOS rather than by newly launched replacements. This shows that the value of OOS lies not only in the financial transfer between stakeholders but also in its ability to reduce dependence on new satellite production and launch. In a broader sense, this may also contribute to more sustainable use of the space environment, although such implications lie beyond the scope of the present model.

Third, among the OOS-related parameters, the serviceable-failure ratio $r_{\oos}$ is the only one that materially changes the constellation operator's preferred maintenance architecture. Increasing $r_{\oos}$ from 0.1 to 0.5 reduces the minimum $\AMC$ from \$873.4M/year to \$661.9M/year and increases the maximum $\AP$ from \$51.7M/year to \$263.1M/year, while also shifting the operator-side decisions farther from the benchmark solution without OOS and toward stronger reliance on OOS. By contrast, variations in $\scost{\oos,\min}$, $\rrate{\oos,\mathrm{ideal}}$, $\alpha_1$, and $\alpha_2$ mainly shift or reshape the Pareto frontier without materially changing the operator-side configuration. From a managerial standpoint, this implies that expanding the set of recoverable failure modes may be at least as important as reducing service cost or improving nominal responsiveness.

\subsubsection{Limitations}
These findings should be interpreted in light of several limitations. First, the cost-responsiveness relationship is imposed exogenously through a specific parametric functional form. Although its qualitative shape is broadly plausible, given the underlying effort/cost trade-off and its consistency with prior work, it may not accurately represent future OOS technologies, operational concepts, or market structures. Second, recovered satellites are assumed to return to an ``as-good-as-new'' state, which likely overstates the long-run benefit of repeated servicing. Third, the analysis adopts a two-stakeholder bi-objective optimization framework and therefore abstracts from bargaining, competition among OOS providers, contracting frictions, and other strategic effects that may determine which solution is realized in practice. As a result, the realized outcome need not coincide with the Pareto-optimal frontier identified by the present model and may, in practice, lie outside it.
}

\section{Conclusions}\label{section 6}
This study introduces an OOS-integrated maintenance strategy for large-scale satellite constellations, in which operators procure spare satellites simultaneously by deploying new units and restoring failed satellites using OOS. A three-echelon supply chain, consisting of in-plane spares (first), parking spares (second), and ground spares (third), is considered. The system forms a closed-loop structure connected to the in-plane spares, where failed satellites are recovered through OOS, and the restored units are used to replenish the in-plane spares.

On this foundation, we develop analytical models for inventory parameters that represent the average behavior of the supply chain over an infinite time horizon, under parametric replenishment policies. Based on these models, we introduce optimization problems relevant to two key stakeholders: the constellation operator and the OOS provider. To capture the OOS provider's decision-making context, we further introduce a functional relationship between OOS cost and responsiveness, characterized by a strictly decreasing curve with a knee point and two asymptotes parallel to the coordinate axes, implying cost and performance bounds. \rev{We then formulate the interaction between the two stakeholders---the constellation operator and the OOS provider---as a bi-objective optimization problem to identify the family of Pareto-optimal solutions. The resulting models and optimization framework are applied to realistic large-scale constellation maintenance scenarios under various OOS-related parameter settings.

The case study yields three main findings. First, within a given scenario, the Pareto-optimal set is generally generated by an almost invariant maintenance strategy on the constellation operator's side, while most variation along the Pareto frontier is driven by the provider-side pricing decision. Second, integrating OOS improves the maintenance system without degrading the required service levels, indicating that its value lies not only in the economic exchange between stakeholders but also in reducing dependence on conventional replacement through new satellite production and launch. Third, among the OOS-related parameters, the fraction of serviceable failure plays the most important structural role, as it is the only parameter that materially changes the constellation operator's preferred maintenance architecture, while the remaining parameters mainly shift or reshape the Pareto frontier without substantially altering the operator-side decisions.

Future work should therefore proceed in three main directions. First, the cost-responsiveness relationship should be refined using more realistic, technology- and mission-specific models, especially as OOS-related technologies, policies, and market conditions continue to evolve. Second, the recovery process could be extended to account for imperfect restoration and state-dependent degradation. Third, the present bi-objective framework should be expanded to incorporate more realistic strategic interactions, such as bargaining and competition among multiple OOS providers.}

\section*{Acknowledgments}
This work was supported by the National Research Foundation of Korea (NRF) grant funded by the Korea government (MSIT) (RS-2025-02213804 \& RS-2025-00563732). The authors thank the NRF for the support of this work.

\rev{

}
\end{document}